\documentclass[aos,nameyear,dvips]{arximspdf}

\doi{10.1214/10-AOS825}
\volume{38}
\issue{6}
\pubyear{2010}
\firstpage{3660}
\lastpage{3695}

\makeatletter

\newproclaim{Example}{Example}
\newproclaim{Remarks}{Remarks}

\newtheorem{theorem}{Theorem}
\newtheorem{proposition}[theorem]{Proposition}
\newtheorem{lemma}[theorem]{Lemma}

\newcommand{\bE}{\mathbb{E}}
\newcommand{\bR}{\mathbb{R}}
\newcommand{\Ecal}{\mathcal E}
\newcommand{\Gcal}{\mathcal G}
\newcommand{\Hcal}{\mathcal H}
\newcommand{\Kcal}{\mathcal K}
\newcommand{\argmin}{\mathop{\arg\min}}

\makeatother

\begin{document}
\begin{frontmatter}

\title{Sparsity in multiple kernel learning}
\runtitle{Sparsity in multiple kernel learning}

\begin{aug}
\author[A]{\fnms{Vladimir} \snm{Koltchinskii}\corref{}\thanksref{t1}\ead[label=e1]{vlad@math.gatech.edu}} and
\author[B]{\fnms{Ming} \snm{Yuan}\thanksref{t2}\ead[label=e2]{myuan@isye.gatech.edu}}
\runauthor{V. Koltchinskii and M. Yuan}
\affiliation{Georgia Institute of Technology}
\address[A]{School of Mathematics\\
Georgia Institute of Technology\\
Atlanta, Georgia 30332-0160\\
USA\\
\printead{e1}} 
\address[B]{School of Industrial\\
\quad and Systems Engineering\\
Georgia Institute of Technology\\
Atlanta, Georgia 30332-0205\\
USA\\
\printead{e2}}
\end{aug}

\thankstext{t1}{Supported in part by
NSF Grants MPSA-MCS-0624841, DMS-09-06880 and CCF-0808863.}

\thankstext{t2}{Supported in part by NSF Grants MPSA-MCS-0624841
and DMS-08-46234.}

\received{\smonth{11} \syear{2009}}
\revised{\smonth{3} \syear{2010}}

%
\begin{abstract}
The problem of multiple kernel learning based on penalized empirical
risk minimization is discussed. The complexity penalty is determined
jointly by the empirical $L_2$ norms and the reproducing kernel Hilbert
space (RKHS) norms induced by the kernels with a data-driven choice of
regularization parameters. The main focus is on the case when the total
number of kernels is large, but only a relatively small number of them
is needed to represent the target function, so that the problem is
sparse. The goal is to establish oracle inequalities for the excess
risk of the resulting prediction rule showing that the method is
adaptive both to the unknown design distribution and to the sparsity of
the problem.
\end{abstract}

%

\end{frontmatter}

\section{Introduction}

Let $(X_i,Y_i), i=1,\ldots,n$ be independent copies of a random couple
$(X,Y)$ with values in $S\times T$,
where $S$ is a measurable space with
$\sigma$-algebra $\mathcal{A}$ (typically, $S$ is
a compact subset of a finite-dimensional Euclidean space)
and $T$ is a Borel subset of $\bR$. In what follows, $P$ will denote the
distribution of $(X,Y)$ and $\Pi$ the distribution of $X$. The corresponding
empirical distributions, based on $(X_1,Y_1),\ldots,(X_n,Y_n)$ and
on $(X_1,\ldots, X_n)$, will be denoted by $P_n$ and $\Pi_n$,
respectively. For a measurable function $g\dvtx S\times T\mapsto\bR$,
we denote
\[
Pg := \int_{S\times T}g\,dP={\bE}g(X,Y)
\quad\mbox{and}\quad
P_n g := \int_{S\times T}g\,dP_n=n^{-1}\sum_{j=1}^n g(X_j,Y_j).
\]
Similarly, we use the notations $\Pi f$ and $\Pi_n f$ for the integrals
of a function $f\dvtx S\mapsto{\bR}$ with respect to the measures $\Pi$ and
$\Pi_n$.

The goal of prediction is to learn ``a reasonably good''
prediction rule $f\dvtx S\to\bR$ from the empirical data
$\{(X_i, Y_i)\dvtx i=1,2,\ldots, n\}$.
To be more specific, consider a loss function
$\ell\dvtx T\times\bR\to\bR_{+}$ and define the risk
of a prediction rule $f$ as
\[
P(\ell\circ f)=\bE\ell(Y,f(X)),
\]
where
$(\ell\circ f)(x,y)=\ell(y,f(x))$.
An optimal
prediction rule with respect to this loss is defined as
\[
f_\ast=\argmin_{f\dvtx S\mapsto\bR} P(\ell\circ f),
\]
where the minimization is taken over all measurable functions
and, for simplicity, it is assumed that the minimum is attained.
\textit{The excess risk} of a prediction rule $f$ is defined as
\[
\mathcal{E}(\ell\circ f):=P(\ell\circ f)-P(\ell\circ f_{\ast}).
\]

Throughout the paper, the notation $a\asymp b$ means
that there exists a numerical constant $c>0$ such that $c^{-1}\leq
\frac
{a}{b}\leq c$.
By ``numerical constants'' we usually mean real numbers whose precise
values are not
necessarily specified, or, sometimes, constants that might depend on
the characteristics of the
problem that are of little interest to us (e.g., some constants
that depend only on the loss function).

\subsection{Learning in reproducing kernel Hilbert spaces}\label{sec11}
Let ${\Hcal}_K$ be a reproducing kernel Hilbert
space (RKHS) associated with a symmetric nonnegatively definite
kernel
$K\dvtx S\times S\to\bR$ such that for any $x\in S$,
$K_x(\cdot):=K(\cdot, x)\in\Hcal_K$
and $f(x)=\langle f, K_x\rangle_{\Hcal_K}$ for all $f\in{\Hcal}_K$
[\citet{Aronszajn1950}]. If it is known that if $f_{\ast}\in
\mathcal
{H}_K$ and
$\|f_{\ast}\|_{\mathcal{H}_K}\leq1$, then it is natural to estimate
$f_{\ast}$ by a solution $\hat f$
of the following empirical risk minimization problem:
%
%
\begin{equation}
\label{simple_ERM}
\hat f := \argmin_{\|f\|_{\mathcal{H}_K}\leq1}
{1\over n}\sum_{i=1}^n\ell(Y_i,f(X_i)).
\end{equation}
The size of the excess risk $\mathcal{E}(\ell\circ\hat f)$
of such an empirical solution
depends on the ``smoothness'' of functions in the RKHS $\mathcal{H}_K$.
A natural notion of ``smoothness'' in this context is related
to the unknown design distribution $\Pi$.
Namely, let $T_{K}$ be the integral operator from $L_2(\Pi)$
into $L_2(\Pi)$ with kernel $K$. Under a standard assumption
that the kernel $K$ is square integrable (in the theory of RKHS
it is usually even assumed that $S$ is compact and $K$ is continuous),
the operator $T_K$ is compact and its spectrum is discrete.
If $\{\lambda_k\}$ is the sequence of the eigenvalues
(arranged in decreasing order) of $T_K$
and $\{\phi_k\}$ is the corresponding $L_2(\Pi)$-orthonormal sequence of
eigenfunctions, then it is well known that the RKHS-norms\vadjust{\goodbreak}
of functions from the linear span of $\{\phi_k\}$ can be written as
\[
\|f\|_{\mathcal{H}_K}^2 = \sum_{k\geq1}
\frac{|\langle f, \phi_k\rangle_{L_2(\Pi)}|^2}{\lambda_k},z
\]
which means that the ``smoothness'' of functions in $\mathcal{H}_K$
depends on the rate of decay of eigenvalues $\lambda_k$ that, in turn,
depends on the design distribution~$\Pi$. It is also clear that the
unit balls in the RKHS ${\mathcal H}_K$ are ellipsoids in the space
$L_2(\Pi)$ with ``axes'' $\sqrt{\lambda_k}$.

It was shown by \citet{Mendelson2002} that the  function
\[
\breve\gamma_n (\delta):=
\biggl(n^{-1}\sum_{k\geq1}(\lambda_k\wedge\delta^2)
\biggr)^{1/2},\qquad
\delta\in[0,1],
\]
provides tight upper and lower bounds (up to constants) on localized
Rade\-macher complexities of the unit ball
in $\mathcal{H}_K$ and
plays an important role in the analysis of the empirical risk\vspace*{1pt} minimization
problem (\ref{simple_ERM}). It is easy to see that the function
$\breve
\gamma_n^2(\sqrt{\delta})$ is concave,
$\breve\gamma_n (0)=0$ and, as a consequence, $\breve\gamma
_n(\delta
)/\delta$ is a
decreasing function of $\delta$ and $\breve\gamma_n(\delta)/\delta^2$
is strictly decreasing.
Hence, there exists unique positive
solution of the equation $\breve\gamma_n(\delta)=\delta^2$. If
$\bar
\delta_n $
denotes this solution, then the results of \citet{Mendelson2002}
imply that
with some constant $C>0$ and with probability at least $1-e^{-t}$
\[
\mathcal{E}(\ell\circ\hat f)\leq C \biggl(\bar\delta_n^2 + \frac
{t}{n} \biggr).
\]
The size of the quantity $\bar\delta_n^2$ involved in this upper bound
on the excess risk depends on the rate of decay of the eigenvalues
$\lambda_k$ as $k\to\infty$. In particular, if $\lambda_k \asymp
k^{-2\beta}$
for some $\beta>1/2$, then it is easy to see that
$
\breve\gamma_n(\delta)\asymp n^{-1/2}\delta^{1-{1/(2\beta)}}
$
and $\bar\delta_n^2 \asymp n^{-2\beta/(2\beta+1)}$.
Recall that unit balls in ${\mathcal H}_K$ are ellipsoids
in $L_2(\Pi)$ with ``axes'' of the order $k^{-\beta}$ and
it is well known that,
in a variety of estimation problems,
$n^{-2\beta/(2\beta+1)}$ represents minimax convergence rates of the
squared $L_2$-risk for functions from such ellipsoids
(e.g., from Sobolev balls of smoothness
$\beta$), as in famous Pinsker's theorem [see, e.g., \citet{Tsybakov09},
Chapter 3].
\begin{Example*}
Sobolev spaces $W^{\alpha,2}(G), G\subset{\mathbb
R}^d$ of smoothness $\alpha>d/2$ is a well-known class of concrete
examples of
RKHS. Let ${\mathbb T}^d, d\geq1$ denote the $d$-dimensional torus
and let $\Pi$ be the uniform distribution in ${\mathbb T}^d$.
It is easy to check that, for all $\alpha>d/2$, the Sobolev space
$W^{\alpha,2}({\mathbb T}^d)$ is
an RKHS generated by the kernel $K(x,y)=k(x-y), x,y\in{\mathbb T}$,
where the function $k\in L_2({\mathbb T}^d)$ is defined by its Fourier
coefficients
\[
\hat k_n=(|n|^2+1)^{-\alpha},\qquad n=(n_1,\ldots, n_d)\in{\mathbb
Z}^d,\qquad
|n|^2:=n_1^2+\cdots+n_d^2.
\]
In this case, the eigenfunctions of the operator $T_K$ are the
functions of the Fourier basis and its eigenvalues are the numbers $\{
(|n|^2+1)^{-\alpha}\dvtx n\in{\mathbb Z}^d\}$.\vadjust{\goodbreak}
For $d=1$ and $\alpha>1/2$, we have
$\lambda_k\asymp k^{-2\alpha}$ (recall that\vspace*{1pt} $\{\lambda_k\}$
are the eigenvalues arranged in decreasing order)
so, $\beta=\alpha$ and
$\bar\delta_n^2\asymp n^{-2\alpha/(2\alpha+1)}$,
which is a minimax nonparametric convergence rate for Sobolev balls
in $W^{\alpha,2}({\mathbb T})$ [see, e.g., \citet{Tsybakov09},
Theorem 2.9].
More generally, for arbitrary $d\geq1$ and $\alpha>d/2$, we get
$\beta=\alpha/d$ and $\bar\delta_n^2 \asymp n^{-2\alpha/(2\alpha+d)}$,
which is also\vspace*{-1pt} a minimax optimal convergence rate in this case.
Suppose now that
the distribution $\Pi$ is uniform in a torus ${\mathbb T}^{d'}\subset
{\mathbb T}^d$ of dimension $d'<d$. We will use the same kernel $K$,
but restrict the RKHS $\mathcal{H}_K$ to the torus ${\mathbb T}^{d'}$
of smaller dimension. Let $d^{\prime\prime}=d-d^{\prime}$. For $n\in
{\mathbb Z^d}$, we
will write $n=(n^{\prime}, n^{\prime\prime})$ with $n^{\prime}\in
{\mathbb Z}^{d'},
n^{\prime\prime}\in{\mathbb Z}^{d^{\prime\prime}}$.
It is easy to prove that the eigenvalues of
the operator $T_K$ become in this case
\[
\sum_{n^{\prime\prime}\in{\mathbb Z}^{d^{\prime\prime}}}(|n^{\prime
}|^2+|n^{\prime\prime}|^2+1)^{-\alpha
}\asymp(|n^{\prime}|^2+1)^{-(\alpha-d^{\prime\prime}/2)}.
\]
Due to this fact,\vspace*{2pt} the norm of the space $\mathcal{H}_K$ (restricted to
${\mathbb T}^{d^{\prime}}$) is equivalent to the norm of the Sobolev
space $W^{\alpha-d^{\prime\prime}/2,2}({\mathbb T}^{d'})$.
Since the eigenvalues of the operator
$T_K$ coincide, up to a constant, with the numbers $\{
(|n^{\prime}|^2+1)^{-(\alpha-d^{\prime\prime}/2)}\dvtx n^{\prime}\in
{\mathbb Z}^{d'}\}$,
we get $\bar\delta_n^2 \asymp n^{-({2\alpha-d^{\prime\prime
}})/({2\alpha
-d^{\prime\prime}+d^{\prime}})}$ [which is again the minimax convergence
rate for Sobolev balls in $W^{\alpha-d^{\prime\prime}/2,2}({\mathbb T}^{d'})$].
In the case of more general
design distributions $\Pi$, the rate of decay of the eigenvalues
$\lambda_k$ and the corresponding size of the excess risk
bound $\bar\delta_n^2$ depends on $\Pi$. If, for instance,
$\Pi$ is supported in a submanifold $S\subset{\mathbb T}^d$
of dimension $\operatorname{dim}(S)<d$, the rate of convergence of
$\bar\delta_n^2$
to $0$ depends on the dimension of the submanifold $S$ rather than
on the dimension of the ambient space ${\mathbb T}^d$.

Using the properties\vspace*{1pt} of the function $\breve\gamma_n$, in particular,
the fact that $\breve\gamma_n(\delta)/\delta$ is decreasing,
it is easy to observe that
$
\breve\gamma_n(\delta)\leq\bar\delta_n \delta+ \bar\delta
_n^2,
\delta\in(0,1].
$
Moreover, if $\breve\epsilon=\breve\epsilon(K)$
denotes
the smallest value of $\epsilon$ such that the linear function
$\epsilon\delta+ \epsilon^2, \delta\in(0,1]$ provides an upper bound
for the
function $\breve\gamma_n(\delta), \delta\in(0,1]$, then
$
\breve\epsilon\leq\bar\delta_n \leq2(\sqrt{5}-1)^{-1}\breve
\epsilon.
$
Note that $\breve\epsilon$ also depends on $n$, but
we do not have to emphasize this dependence in
the notation since, in what follows, $n$ is fixed.
Based on the observations above, the quantity $\bar\delta_n$ coincides
(up to a numerical constant) with the
slope $\breve\epsilon$ of the ``smallest linear majorant'' of the form
$\epsilon\delta+\epsilon^2$
of the function $\breve\gamma_n(\delta)$.
This interpretation of $\bar\delta_n$ is of some
importance in the design of complexity penalties used in this paper.
\end{Example*}

\subsection{Sparse recovery via regularization}
Instead of minimizing the empirical risk over an RKHS-ball [as in problem
(\ref{simple_ERM})], it is very common to define the estimator $\hat f$
of the target function $f_{\ast}$ as a solution of the penalized empirical
risk minimization problem of the form
%
%
\begin{equation}
\label{AA}
\hat f:=\argmin_{f\in\Hcal} \Biggl[{1\over n}\sum_{i=1}^n\ell
(Y_i,f(X_i))+
\epsilon\|f\|_{{\Hcal}_K}^\alpha\Biggr],
\end{equation}
where $\epsilon>0$ is a tuning parameter that
balances the tradeoff between the empirical risk and the ``smoothness''
of the estimate and, most often, $\alpha=2$ (sometimes, \mbox{$\alpha=1$}).
The properties of the estimator $\hat f$ has been studied extensively.
In particular, it was possible to derive probabilistic bounds on the
excess risk $\mathcal{E}(\ell\circ\hat f)$ (oracle inequalities) with
the control
of the random error in terms of the rate of decay of the eigenvalues
$\{\lambda_k\}$, or, equivalently, in terms of the function $\breve
\gamma_n$
[see, e.g., \citet{BlanchardBousquetMassart2008}].\looseness=-1

In the recent years, there has been a lot of interest in a data dependent
choice of kernel $K$ in this type of problems. In particular, given
a finite (possibly large) dictionary
$\{K_j\dvtx j=1,2,\ldots, N\}$ of symmetric nonnegatively definite kernels
on~$S$,
one can try to find a ``good'' kernel $K$ as a convex combination
of the kernels from the dictionary:
%
%
\begin{equation}
K\in\Kcal:= \Biggl\{\sum_{j=1}^N\theta_jK_j\dvtx
\theta_j\ge0, \theta_1+\cdots+\theta_N=1 \Biggr\}.
\end{equation}
The coefficients of $K$ need to be estimated
from the training data along with the prediction rule.
Using this approach for problem (\ref{AA}) with $\alpha=1$ leads to the
following optimization problem:
%
%
\begin{equation}
\label{AB}
\hat f:=\mathop{\mathop{\arg\min}_{f\in\mathcal{H}_K}}_{K\in\mathcal{K}}
\bigl(P_n(\ell\circ f)+\epsilon\|f\|_{{\Hcal}_K} \bigr).
\end{equation}
This learning problem, often referred
to as the multiple kernel learning,
has been studied recently by \citet{BousquetHerrmann2003},
\citet{CrammerKeshetSinger2003},
\citet{Lanckrietetal04}, \citet{MicchelliPontil2005},
\citet{LinZhang2006},
\citet{SBd06}, \citet{Bach2008} and \citet
{KoltchinskiiYuan2008} among others.
In particular [see, e.g., \citet{MicchelliPontil2005}], problem
(\ref{AB})
is equivalent to the following:
%
%
\begin{eqnarray}
\label{LASSO_inf}
(\hat f_1,\ldots, \hat f_N)&:=&\argmin_{f_j\in\mathcal{H}_{K_j},
j=1,\ldots, N} \Biggl(P_n\bigl(\ell\circ(f_1+\cdots+f_N)\bigr)\nonumber\\[-8pt]\\[-8pt]
&&\hspace*{104.5pt}{}+\epsilon
\sum_{j=1}^N \|f_j\|_{\mathcal{H}_{K_j}} \Biggr),\nonumber
\end{eqnarray}
which is an infinite-dimensional version of LASSO-type penalization.
\citet{KoltchinskiiYuan2008} studied this method in the case when the
dictionary is large, but the target function $f_{\ast}$ has a
``sparse representation'' in terms of a relatively small subset
of kernels $\{K_j\dvtx j\in J\}$. It was shown that this method is adaptive
to sparsity extending well-known properties of LASSO to this
infinite-dimensional framework.

In this paper, we study a different approach to the multiple kernel learning.
It is closer to the recent work on ``\textit{sparse additive models}''
[see, e.g., \citet{Ravikumaretal08} and\vadjust{\goodbreak}
\citet{MGB09}]
and it is based on a ``double penalization''
with a combination of empirical $L_2$-norms (used
to enforce the sparsity of the solution) and RKHS-norms (used to enforce
the ``smoothness'' of the components). Moreover, we suggest a data-driven
method of choosing the values of regularization parameters that is adaptive
to unknown smoothness of the components (determined by the behavior of
distribution dependent eigenvalues of the kernels).

Let $\Hcal_j:=\mathcal{H}_{K_j}, j=1,\ldots, N$. Denote
$\mathcal{H}:= \mathrm{l.s.} (\bigcup_{j=1}^N \mathcal{H}_j )$
(``l.s.'' meaning ``the linear span''), and
\[
\mathcal{H}^{(N)}:= \{(h_1,\ldots, h_N)\dvtx h_j\in\mathcal{H}_j,
j=1,\ldots, N \}.
\]
Note that $f\in\mathcal{H}$ if and only if there exists an
additive\vspace*{1pt}
representation (possibly, nonunique)
$
f=f_1+\cdots+f_N,
$
where $f_j\in\mathcal{H}_j$, $j=1,\ldots, N$.
Also, $\mathcal{H}^{(N)}$ has a natural structure of a linear space
and it can be equipped with the following inner product:
\[
\langle(f_1,\ldots, f_N), (g_1,\ldots, g_N)\rangle_{\mathcal{H}^{(N)}}
:=\sum_{j=1}^N \langle f_j, g_j\rangle_{\mathcal{H}_j}
\]
to become \textit{the direct sum} of Hilbert spaces $\mathcal{H}_j,
j=1,\ldots, N$.

Given a convex subset $D\subset\mathcal{H}^{(N)}$,
consider the following penalized empirical risk minimization problem:
%
%
\begin{eqnarray}
\label{main_ERM}
(\hat{f}_1,\ldots, \hat{f}_N )&=&\argmin_{(f_1,\ldots, f_N)
\in D}
\Biggl[P_n\bigl(\ell\circ(f_1+\cdots+f_N)\bigr)\nonumber\\[-8pt]\\[-8pt]
&&\hspace*{54.4pt}{}+\sum_{j=1}^N \bigl({\epsilon
}_j\|
f_j\|_{L_2(\Pi_n)}
+{\epsilon}_j^2\|f_j\|_{\Hcal_j} \bigr) \Biggr].\nonumber
\end{eqnarray}
Note that for special choices of set $D$, for instance, for
$D:=\{(f_1,\ldots, f_N)\dvtx f_j\in\mathcal{H}_j, \|f_j\|_{\mathcal
{H}_j}\leq R_j\}$
for some $R_j>0, j=1,\ldots, N$, one can replace each component $f_j$
involved in the optimization problem by its orthogonal projections
in $\mathcal{H}_j$ onto the linear span of the functions $\{K_j(\cdot,
X_i), i=1,\ldots, n\}$ and reduce the problem to a convex optimization over
a finite-dimensional space (of dimension $n N$).

The complexity penalty in the problem (\ref{main_ERM}) is based on two
norms of the components $f_j$ of an additive representation:
the empirical $L_2$-norm, $\|f_j\|_{L_2(\Pi_n)}$,
with regularization parameter $\epsilon_j$,
and an RKHS-norm, $\|f_j\|_{\mathcal{H}_j}$, with regularization
parameter $\epsilon_j^2$. The empirical $L_2$-norm (the lighter norm)
is used to enforce the sparsity of the solution whereas the RKHS norms
(the heavier norms) are used to enforce the ``smoothness'' of the
components. This is similar to the approach taken in
\citet{MGB09} in the
context of classical additive models, that is, in the case
when $S:=[0,1]^N$, $\mathcal{H}_j:=W^{\alpha,2}([0,1])$ for some
smoothness $\alpha>1/2$ and the space $\mathcal{H}_j$ is a space of\vadjust{\goodbreak}
functions depending on the $j$th variable. In this case, the regularization
parameters $\epsilon_j$ are equal (up to a constant) to $n^{-\alpha
/(2\alpha+1)}$. The quantity $\epsilon_j^2$, used in the
``smoothness part'' of the penalty, coincides with
the minimax convergence rate in a one component smooth problem.
At the same time, the quantity $\epsilon_j$, used
in the ``sparsity part'' of the penalty, is equal to the
square root of the minimax rate (which is similar to the choice of
regularization parameter in standard sparse recovery methods such as
LASSO). This choice of regularization parameters results in
the excess risk of the order $dn^{-2\alpha/(2\alpha+1)}$, where
$d$ is the number of components of the target function (the
degree of sparsity of the problem).

The framework of multiple kernel learning considered in this
paper includes many generalized versions of classical additive
models. For instance, one can think of the case when
$S:=[0,1]^{m_1}\times\cdots\times[0,1]^{m_N}$ and
$\mathcal{H}_j=W^{\alpha,2}([0,1]^{m_j})$ is a
space of functions depending on the $j$th block of variables.
In this case, a~proper choice of regularization parameters
(for uniform design distribution)
would be $\epsilon_j=n^{-\alpha/(2\alpha+m_j)}, j=1,\ldots, N$
(so, these parameters and the error rates for different components of
the model are different). It should be also clear from the discussion
in Section~\ref{sec11} that, if the design distribution $\Pi$ is unknown,
the minimax convergence rates for the one component problems are also
unknown. For instance, if the projections of design points on the
cubes $[0,1]^{m_j}$ are distributed in lower-dimensional submanifolds
of these cubes, then the unknown dimensions of the submanifolds rather
than the dimensions $m_j$ would be involved in the minimax rates and
in the regularization parameters $\epsilon_j$. Because of this, data
driven choice of regularization parameters $\epsilon_j$ that provides
adaptation to the unknown design distribution $\Pi$ and to the unknown
``smoothness'' of the components (related to this distribution)
is a major issue in multiple kernel learning. From this point of view,
even in the case of classical additive models, the choice of
regularization parameters that is based only on Sobolev type smoothness
and ignores the design distribution is not adaptive.
Note that, in the infinite-dimensional LASSO studied in \citet
{KoltchinskiiYuan2008}, the
regularization parameter $\epsilon$ is chosen the same way as
in the classical LASSO ($\epsilon\asymp\sqrt{\frac{\log N}{n}}$),
so, it is not related to the smoothness of the components. However,
the oracle inequalities proved in \citet{KoltchinskiiYuan2008}
give correct size of the excess risk only for special choices
of kernels that depend on unknown ``smoothness'' of the components
of the target function~$f_{\ast}$, so, this method is not adaptive either.


\subsection{Adaptive choice of regularization parameters}
Denote
\[
\hat K_j := \biggl(\frac{K_j(X_l, X_k)}{n} \biggr)_{l,k=1,n}.\vadjust{\goodbreak}
\]
This $n\times n$ Gram matrix can be viewed as an empirical version
of the integral operator $T_{K_j}$
from $L_2(\Pi)$ into $L_2(\Pi)$ with kernel $K_j$.
Denote\vadjust{\goodbreak} $\hat\lambda_k^{(j)}, k=1,2,\ldots,$
the eigenvalues of $\hat K_j$ arranged in decreasing order.
We also use the notation $\lambda_k^{(j)}, k=1,2,\ldots,$
for the eigenvalues of the operator $T_{K_j}\dvtx L_2(\Pi)\mapsto
L_2(\Pi)$
with kernel $K_j$ arranged in decreasing order.
Define functions $\breve\gamma_n^{(j)}, \hat\gamma_n^{(j)}$,
\[
\breve\gamma_n^{(j)}(\delta):=
\Biggl({1\over n}\sum_{k=1}^n
\bigl(\lambda_k^{(j)}\wedge\delta^2\bigr) \Biggr)^{1/2}
\quad\mbox{and}\quad
\hat\gamma_n^{(j)}(\delta):=
\Biggl({1\over n}\sum_{k=1}^n
\bigl(\hat\lambda_k^{(j)}\wedge\delta^2\bigr) \Biggr)^{1/2},
\]
and, for a fixed given $A\geq1$, let
%
%
\begin{equation}
\label{regul_par}
\hat{\epsilon}_j
:=\inf\Biggl\{\epsilon\ge\sqrt{A\log N\over n}\dvtx
\hat\gamma_n^{(j)}(\delta)\leq
\epsilon\delta+ \epsilon^2, \forall\delta\in(0,1]
\Biggr\}.
\end{equation}
One can view $\hat\epsilon_j$ as an empirical estimate of the quantity
$\breve\epsilon_j=\breve\epsilon(K_j)$ that
(as we have already pointed out)
plays a crucial role in the bounds on the excess risk in empirical
risk minimization problems in the RKHS context.
In fact, since most often $\breve\epsilon_j\geq\sqrt{{A\log N}/{n}}$,
we will redefine this quantity as
%
%
\begin{equation}
\label{regul_par_01}
\breve{\epsilon}_j
:=\inf\Biggl\{\epsilon\ge\sqrt{A\log N\over n}\dvtx
\breve\gamma_n^{(j)}(\delta)\leq
\epsilon\delta+ \epsilon^2, \forall\delta\in(0,1]
\Biggr\}.
\end{equation}

We will use the following values of regularization parameters in problem~(\ref{main_ERM}):
$\epsilon_j=\tau\hat{\epsilon}_j$,
where $\tau$ is a sufficiently large constant.

It should be emphasized that the structure of complexity penalty
and the choice of regularization parameters in (\ref{main_ERM})
are closely related to the following bound on Rademacher processes
indexed by functions from an RKHS $\mathcal{H}_K$: with a high probability,
for all $h\in\mathcal{H}_K$,
\[
|R_n(h)|\leq C \bigl[\breve\epsilon(K)\|h\|_{L_2(\Pi)}+
\breve\epsilon^2(K)\|h\|_{\mathcal{H}_K} \bigr].
\]
Such bounds follow from the results of Section~\ref{sec:sam2pop} and
they provide
a way to prove sparsity oracle inequalities for the estimators
(\ref{main_ERM}).
The Rademacher process is defined as
\[
R_n(f):=n^{-1}\sum_{j=1}^n \varepsilon_j f(X_j),
\]
where $\{\varepsilon_j\}$ is a sequence of
i.i.d. Rademacher random variables (taking values $+1$ and $-1$ with
probability $1/2$ each) independent of $\{X_j\}$.

We will use several basic facts of the empirical processes theory
throughout the paper. They include symmetrization inequalities and contraction
(comparison) inequalities for Rademacher processes that can be found in the
books of \citet{LedouxTalagrand1991} and \citet{VW96}.
We also use Talagrand's concentration inequality for empirical processes
[see, \citet{Talagrand1996}, \citet{Bousquet2002}].\vadjust{\goodbreak}

The main goal of the paper is to establish oracle
inequalities for the excess risk of the estimator
$\hat{f}=\hat{f}_1+\cdots+\hat{f}_N$.
In these inequalities, the excess risk of $\hat f$ is compared
with the excess risk of an oracle
$f:=f_1+\cdots+f_N, (f_1,\ldots, f_N)\in D$
with an error term depending on the degree of sparsity
of the oracle, that is, on the number of nonzero components
$f_j\in\mathcal{H}_j$ in its additive representation.
The oracle inequalities will be stated in the next section.
Their proof relies on probabilistic
bounds for empirical $L_2$-norms
and data dependent regularization parameters $\hat{\epsilon}_j$.
The results of Section~\ref{sec:sam2pop}
show that they
can be bounded by their respective population counterparts.
Using these tools and some bounds on empirical processes derived
in Section~\ref{sec:empirical_process},
we prove in Section~\ref{sec:mainres} the oracle inequalities
for the estimator~$\hat{f}$.

\section{Oracle inequalities}
\label{sec:oracle}

Considering the problem\vspace*{1pt} in the case when
the domain $D$ of (\ref{main_ERM})
is not bounded, say, $D=\Hcal^{(N)}$,
leads to additional technical complications and
might require some changes in the estimation
procedure. To avoid this, we assume below that $D$ is a bounded convex
subset of $\mathcal{H}^{(N)}$.
It will be also assumed that, for all $j=1,\ldots, N$,
$
\sup_{x\in S}K_j(x,x)\leq1,
$
which, by elementary properties of RKHS, implies that
$
\|f_j\|_{L_{\infty}}\leq\|f_j\|_{\mathcal{H}_j}, j=1,\ldots, N.
$
Because of this,
\[
R_D:=\sup_{(f_1,\ldots, f_N)\in D}\|f_1+\cdots+ f_N\|_{L_{\infty
}}<+\infty.
\]
Denote
$
R_D^{\ast}:=R_D\vee\|f_{\ast}\|_{L_{\infty}}.
$
We will allow the constants involved in the oracle inequalities
stated and proved below to depend on the value of $R_D^{\ast}$
(so, implicitly,
it is assumed that this value is not too large).

We shall also assume that $N$ is large enough, say, so that
$\log N\ge2\log\log n$.
This assumption is not essential to our
development and is in place to avoid
an extra term of the order $n^{-1}\log\log n$ in our risk bounds.

\subsection{Loss functions of quadratic type}
We will formulate the \textit{assumptions on the loss function $\ell$}.
The main assumption is that, for all $y\in T$,
$\ell(y,\cdot)$ is a nonnegative convex function. In addition, we will
assume that $\ell(y,0), y\in T$ is uniformly bounded from above by a
numerical constant. Moreover, suppose that, for all $y\in T$,
$\ell(y,\cdot)$ is twice continuously differentiable and its
first and second derivatives are uniformly bounded in $T\times
[-R_D^{\ast}, R_D^{\ast}]$.
Denote
%
%
\begin{equation}
\label{SC}
m(R):=\frac{1}{2}\inf_{y\in T}
\inf_{|u|\leq R}{\partial^2 \ell(y,u)\over\partial u^2},\qquad
M(R):=\frac{1}{2}\sup_{y\in T}
\sup_{|u|\leq R}{\partial^2 \ell(y,u)\over\partial u^2}
\end{equation}
and let
$
m_{\ast}:=m(R_D^{\ast}), M_{\ast}:=M(R_D^{\ast}).
$
We will assume that $m_{\ast}>0$.

Denote
\[
L_{\ast}:=\sup_{|u|\leq R_D^{\ast}, y\in T} \biggl|{\partial\ell
\over
\partial
u}(y,u) \biggr|.
\]
Clearly, for all $y\in T$, the function $\ell(y,\cdot)$ satisfies
Lipschitz condition with constant $L_{\ast}$.\vadjust{\goodbreak}

The constants $m_{\ast}, M_{\ast}, L_{\ast}$ will appear
in a number of places in what follows.
Without loss of generality, we can also assume that $m_{\ast}\leq1$
and $L_{\ast}\geq1$
(otherwise, $m_{\ast}$ and $L_{\ast}$
can be replaced by a lower bound and an upper bound, resp.).

The loss functions satisfying the assumptions stated above will
be called \textit{the losses of quadratic type}.

If $\ell$ is a loss of quadratic type and $f=f_1+\cdots+f_N,
(f_1,\ldots, f_N)\in D$, then
%
%
\begin{equation}
\label{excess_L_2}
m_{\ast}
\|f-f_{\ast}\|_{L_2(\Pi)}^2\le\mathcal{E}(\ell\circ f) \leq
M_{\ast}\|f-f_{\ast}\|_{L_2(\Pi)}^2.
\end{equation}
This bound easily follows from a simple argument based on
Taylor expansion and it will be used later in the paper.
If
$\mathcal{H}$ is dense in $L_2(\Pi)$, then (\ref{excess_L_2})
implies that
%
%
\begin{equation}
\inf_{f\in\mathcal{H}}P(\ell\circ f)=
\inf_{f\in L_2(\Pi)}P(\ell\circ f)
=P(\ell\circ f_{\ast}).
\end{equation}

The quadratic loss $\ell(y,u):=(y-u)^2$
in the case when $T\subset\bR$ is a bounded set is
one of the main examples of such loss functions.
In this case, $m(R)=1$ for all $R>0$.
In regression problems with a bounded response variable,
more general loss functions of the form $\ell(y,u):=\phi(y-u)$
can be also used,
where $\phi$ is an even nonnegative convex twice
continuously differentiable function with
$\phi^{\prime\prime}$ uniformly bounded
in ${\mathbb R}$, $\phi(0)=0$ and
$\phi^{\prime\prime}(u)>0, u\in{\mathbb R}$.
In classification problems,
the loss functions of the form $\ell(y,u)=\phi(yu)$
are commonly used, with $\phi$ being a
nonnegative decreasing convex twice
continuously differentiable function such that,
again, $\phi^{\prime\prime}$
is uniformly bounded
in ${\mathbb R}$ and $\phi^{\prime\prime}(u)>0, u\in{\mathbb R}$.
The loss function $\phi(u)=\log_2(1+e^{-u})$
(often referred to as the logit loss) is a specific example.

\subsection{Geometry of the dictionary}
Now we introduce several important geometric characteristics
of dictionaries consisting of kernels (or, equivalently,
of RKHS). These characteristics are related to the degree
of ``dependence'' of spaces of random variables $\mathcal{H}_j\subset
L_2(\Pi),
j=1,\ldots, N$ and they will be involved in the oracle
inequalities for the excess risk $\mathcal{E}(\ell\circ\hat f)$.

First, for $J\subset\{1,\ldots, N\}$ and
$b\in[0,+\infty]$,
denote
\[
C_J^{(b)} := \biggl\{(h_1,\ldots, h_N)\in\mathcal{H}^{(N)}\dvtx
\sum_{j\notin J}\|h_j\|_{L_2(\Pi)}\leq b\sum_{j\in J}\|h_j\|
_{L_2(\Pi)}
\biggr\}.
\]
Clearly, the set $C_J^{(b)}$ is a cone in the space $\mathcal{H}^{(N)}$
that consists of vectors $(h_1,\ldots, h_N)$ whose components corresponding
to $j\in J$ ``dominate'' the rest of the components. This family of cones
increases as $b$ increases. For $b=0$, $C_{J}^{(b)}$
coincides with the linear subspace of vectors for which $h_j=0, j\notin J$.
For $b=+\infty$, $C_J^{(b)}$ is the whole space $\mathcal{H}^{(N)}$.

The following quantity will play the most important role:
\begin{eqnarray*}
\beta_{2,b}(J;\Pi):\!&=&\beta_{2,b}(J)\\[-2pt]
:\!&=&
\inf\Biggl\{\beta>0\dvtx
\biggl(\sum_{j\in J}\|h_j\|_{L_2(\Pi)}^2 \biggr)^{1/2}
\leq\beta\Biggl\|\sum_{j=1}^N h_j \Biggr\|_{L_2(\Pi)},\\[-2pt]
&&\hspace*{18.17pt}\hspace*{117pt}
(h_1,\ldots, h_N)\in C_{J}^{(b)}
\Biggr\}.
\end{eqnarray*}
Clearly, $\beta_{2,b}(J;\Pi)$ is a nondecreasing
function of $b$. In the case of ``simple dictionary'' that
consists of one-dimensional spaces similar quantities have
been used in the literature on sparse recovery [see, e.g., Koltchinskii
(\citeyear{Koltchinskii2008},
\citeyear{Koltchinskii09a}, \citeyear{Koltchinskii09b}, \citeyear
{Koltchinskii09c}); Bickel, Ritov and Tsybakov (\citeyear{BickelRitovTsybakov2009})].

The quantity $\beta_{2,b}(J;\Pi)$ can be upper bounded
in terms of some other geometric characteristics that describe how
``dependent'' the spaces of random variables $\mathcal{H}_j\subset
L_2(\Pi)$
are. These characteristics will be introduced below.

Given $h_j\in\mathcal{H}_j, j=1,\ldots, N$,
denote by $\kappa(\{h_j\dvtx j\in J\})$
the minimal eigenvalue of the
Gram matrix $(\langle h_j,h_k\rangle_{L_2(\Pi)})_{j,k\in J}$. Let
%
%
\begin{equation}
\kappa(J):=
\inf\bigl\{\kappa(\{h_j\dvtx j\in J\})\dvtx
h_j\in\mathcal{H}_j, \|h_j\|_{L_2(\Pi)}=1 \bigr\}.
\end{equation}
We will also use the notation
%
%
\begin{equation}
\Hcal_J=\mathrm{l.s.} \biggl(\bigcup_{j\in J}\Hcal_j \biggr).
\end{equation}
The following quantity is the maximal cosine of the angle in the space
$L_2(\Pi)$ between
the vectors in the subspaces $\mathcal{H}_I$ and $\mathcal{H}_J$ for some
$I,J\subset\{1,\ldots, N\}$:
%
%
\begin{equation}\quad
\rho(I,J):=
\sup\biggl\{ \frac{\langle f,g\rangle_{L_2(\Pi)}} {\|f\|_{L_2(\Pi)}
\|g\|_{L_2(\Pi)}}\dvtx f\in\Hcal_I, g\in\Hcal_{J}, f\neq0, g\neq0
\biggr\}.
\end{equation}
Denote
$
\rho(J):=\rho(J,J^c).
$
The quantities $\rho(I,J)$ and $\rho(J)$ are very similar to the
notion of \textit{canonical correlation} in the multivariate statistical
analysis.

There are other important geometric characteristics, frequently used
in the theory of sparse recovery, including so called
``\textit{restricted isometry constants}'' by Candes and Tao (\citeyear{CT07}).
Define $\delta_d (\Pi)$ to be the smallest
$\delta>0$ such that for all
$(h_1,\ldots, h_N)\in\mathcal{H}^{(N)}$ and all $J\subset\{1,\ldots,
N\}$
with $\operatorname{card}(J)=d$,
\[
(1-\delta) \biggl(\sum_{j\in J} \|h_j\|_{L_2(\Pi)}^2 \biggr)^{1/2}
\leq
\biggl\|\sum_{j\in J} h_j \biggr\|_{L_2(\Pi)}
\leq
(1+\delta) \biggl(\sum_{j\in J} \|h_j\|_{L_2(\Pi)}^2 \biggr)^{1/2}.
\]
This condition with a sufficiently small value of $\delta_d(\Pi)$
means that for all choices of $J$ with $\operatorname{card}(J)=d$
the functions in the spaces $\mathcal{H}_j, j\in J$ are ``almost orthogonal''
in $L_2(\Pi)$.

The following simple proposition easily follows from some statements
in Koltchinskii (\citeyear{Koltchinskii09a}, \citeyear
{Koltchinskii09b}), (\citeyear{Koltchinskii2008})
(where the case of simple dictionaries
consisting of one-dimensional spaces $\mathcal{H}_j$ was considered).\vadjust{\goodbreak}
\begin{proposition}
For all $J\subset\{1,\ldots, N\}$,
\[
\beta_{2,\infty}(J;\Pi)\leq\frac{1}{\sqrt{\kappa(J)(1-\rho^2(J))}}.
\]
Also, if $\operatorname{card}(J)=d$ and
$\delta_{3d}(\Pi)\leq\frac{1}{8b}$, then $\beta_{2,b}(J;\Pi)\leq4$.
\end{proposition}

Thus, such quantities as $\beta_{2,\infty}(J;\Pi)$ or $\beta
_{2,b}(J;\Pi)$,
for finite values of $b$, are reasonably small provided
that the spaces of random variables $\mathcal{H}_j, j=1,\ldots, N$
satisfy proper conditions of ``weakness of correlations.''

\subsection{Excess risk bounds}
We are now in a position to formulate our main theorems that provide
oracle inequalities for the excess risk $\mathcal{E}(\ell\circ\hat f)$.
In these theorems, $\mathcal{E}(\ell\circ\hat f)$ will be compared
with the excess risk $\mathcal{E}(\ell\circ f)$ of an oracle
$(f_1,\ldots, f_N)\in D$. Here and in what follows,
$f:=f_1+\cdots+f_N\in\mathcal{H}$.
This is a little abuse of notation: we are ignoring the fact that
such an additive representation of a function $f\in\mathcal{H}$ is not
necessarily unique. In some sense, $f$ denotes both the vector
$(f_1,\ldots, f_N)\in\mathcal{H}^{(N)}$ and the function
$f_1+\cdots+f_N\in\mathcal{H}$.
However, this is not going to cause a confusion
in what follows. We will also use the following notation:
\[
J_f:=\{1\leq j\leq N\dvtx f_j\neq0\} \quad\mbox{and}\quad d(f):=\operatorname{card}(J_f).
\]

The error terms of the oracle inequalities will depend on the quantities
$\breve\epsilon_j=\breve\epsilon(K_j)$ related to the ``smoothness''
properties of the RKHS and also on the geometric
characteristics of the dictionary introduced above.
In the first theorem, we will use the quantity $\beta_{2,\infty
}(J_f;\Pi)$
to characterize the properties of the dictionary. In this case, there will
be no assumptions on the quantities $\breve\epsilon_j$: these quantities
could be of different order for different kernel machines, so, different
components of the additive representation could have different ``smoothness.''
In the second
theorem, we will use a smaller quantity $\beta_{2,b}(J;\Pi)$ for a proper
choice of parameter $b<\infty$. In this case, we will have to make an
additional assumption that $\breve\epsilon_j, j=1,\ldots, N$ are all
of the same order (up to a constant).

In both cases, we
consider penalized empirical risk minimization problem~(\ref{main_ERM}) with data-dependent regularization parameters
$\epsilon_j=\tau\hat\epsilon_j$, where
$\hat\epsilon_j, j=1,\ldots, N$ are defined by (\ref{regul_par}) with
some $A\geq4$ and $\tau\geq B L_{\ast}$ for a numerical constant $B$.
\begin{theorem}
\label{co:homobd}
There exist numerical constants $C_1, C_2>0$ such that, for all
all oracles $(f_1,\ldots, f_N)\in D$,
with probability at least $1-3 N^{-A/2}$,
%
%
\begin{eqnarray}
&&
\Ecal(\ell\circ\hat{f})
+C_1 \Biggl(\tau\sum_{j=1}^N \breve\epsilon_j\|\hat f_j-f_j\|
_{L_2(\Pi)}
+\tau^2\sum_{j=1}^N \breve\epsilon_j^2 \|\hat f_j\|_{\mathcal{H}_j}\Biggr)
\nonumber\\[-8pt]\\[-8pt]
&&\qquad
\le
2\Ecal(\ell\circ f)+C_2\tau^2
\sum_{j\in J_f}\breve{\epsilon}_j^2
\biggl(\frac{\beta_{2,\infty}^2(J_f,\Pi)}{m_{\ast}}+
\|f_j\|_{\Hcal_j} \biggr).\nonumber\vadjust{\goodbreak}
\end{eqnarray}
\end{theorem}

This result means that if there exists an oracle $(f_1,\ldots, f_N)\in D$
such that:
\begin{enumerate}[(a)]
\item[(a)] the excess risk $\mathcal{E}(\ell\circ f)$ is small;
\item[(b)] the spaces $\mathcal{H}_j, j\in J_f$ are not strongly correlated
with the spaces $\mathcal{H}_j, j\notin J_f;$
\item[(c)] $\mathcal{H}_j, j\in J_f$
are ``well posed'' in the sense that $\kappa(J_f)$ is not too small;
\item[(d)] $\|f_j\|_{\mathcal{H}_j}, j\in J_f$ are all bounded by a reasonable
constant,
\end{enumerate}
then the excess risk $\mathcal{E}(\ell\circ\hat f)$ is essentially
controlled by $\sum_{j\in J_f}\breve\epsilon_j^2$. At the same time,
the oracle inequality provides a bound on the $L_2(\Pi)$-distances between
the estimated components $\hat f_j$ and the components of the oracle
(of course, everything is under the assumption that the loss is of
quadratic type and $m_{\ast}$ is bounded away from $0$).

Not also that the constant $2$ in front of the excess risk of the oracle
$\Ecal(\ell\circ f)$ can be replaced by $1+\delta$ for any $\delta>0$
with minor modifications of the proof (in this case, the constant $C_2$
depends on $\delta$ and is of the order $1/{\delta}$).

Suppose now that there exists $\breve\epsilon>0$ and a constant
$\Lambda>0$
such that
\[
\Lambda^{-1}\leq\frac{\breve\epsilon_j}{\breve\epsilon}\leq
\Lambda,\qquad
j=1,\ldots, N.
\]

\begin{theorem}
\label{co:homobd-1}
There exist numerical constants $C_1,C_2, b>0$
such that, for all oracles $(f_1,\ldots, f_N)\in D$,
with probability at least $1-3 N^{-A/2}$,
%
%
\begin{eqnarray}
&&
\Ecal(\ell\circ\hat{f})
+\frac{C_1}{\Lambda} \Biggl(\tau\breve\epsilon\sum_{j=1}^N \|\hat
f_j-f_j\|_{L_2(\Pi)}
+\tau^2 \breve\epsilon^2\sum_{j=1}^N \|\hat
f_j\|_{\mathcal{H}_j} \Biggr)
\nonumber\\[-8pt]\\[-8pt]
&&\qquad
\le2\Ecal(\ell\circ f)+
C_2 \Lambda\tau^2\breve
\epsilon^2 \biggl(\frac{\beta_{2,b\Lambda^2}^2(J_f,\Pi)}{m_{\ast
}}d(f) +
\sum_{j\in J_f}\|f_j\|_{\Hcal_j} \biggr).\nonumber
\end{eqnarray}
\end{theorem}

As before,\vspace*{1pt} the constant $2$ in the upper bound can be replaced
by $1+\delta$, but, in this case, the constants $C_2$ and $b$ would be
of the order $\frac{1}{\delta}$.
The meaning of this result is that if
there exists an oracle $(f_1,\ldots, f_N)\in D$
such that:
\begin{enumerate}[(a)]
\item[(a)] the excess risk $\mathcal{E}(\ell\circ f)$ is small;
\item[(b)] the ``restricted isometry'' constant $\delta_{3d}(\Pi)$
is small
for $d=d(f);$
\item[(c)] $\|f_j\|_{\mathcal{H}_j}, j\in J_f$ are all bounded by a reasonable
constant,
\end{enumerate}
then the excess risk $\mathcal{E}(\ell\circ\hat f)$ is essentially
controlled by $d(f)\breve\epsilon^2$. At the same time, the distance
$\sum_{j=1}^N \|\hat f_j-f_j\|_{L_2(\Pi)}$ between the estimator
and the oracle is controlled by $d(f)\breve\epsilon$.
In particular, this implies
that the empirical solution $(\hat f_1,\ldots, \hat f_N)$ is
``approximately sparse'' in the sense that
$\sum_{j\notin J_f}\|\hat f\|_{L_2(\Pi)}$ is of the order
$d(f)\breve
\epsilon$.
\begin{Remarks*} 1. It is easy to check that Theorems~\ref{co:homobd}
and~\ref{co:homobd-1} hold also if one replaces $N$ in the definitions
(\ref{regul_par}) of $\hat\epsilon_j$ and (\ref{regul_par_01}) of
$\breve\epsilon_j$ by an arbitrary $\bar N\geq N$ such that
$\log\bar N\ge2\log\log n $ (a similar condition on $N$ introduced
early in Section~\ref{sec:oracle} is not needed here).
In this case, the probability bounds in the theorems become $1-3{\bar
N}^{-A/2}$.
This change might be of interest if one uses the results for a dictionary
consisting of just one RKHS ($N=1$), which is not the focus of this paper.

2. If the distribution dependent quantities $\breve\epsilon_j,
j=1,\ldots, N$ are known and used as regularization parameters
in (\ref{main_ERM}), the oracle inequalities of Theorems~\ref{co:homobd}
and~\ref{co:homobd-1} also hold (with obvious simplifications of
their proofs). For instance, in the case when $S=[0,1]^N$,
the design distribution $\Pi$ is uniform and, for each $j=1,\ldots, N$,
${\mathcal H}_j$~is a Sobolev space of functions
of smoothness $\alpha>1/2$ depending only on the $j$th variable,
we have $\breve\epsilon_j\asymp n^{-\alpha/(2\alpha+1)}$.
Taking in this case
\[
\epsilon_j=\tau\Biggl(n^{-\alpha/(2\alpha+1)}\vee\sqrt{\frac
{A\log
N}{n}} \Biggr)
\]
would lead to oracle
inequalities for sparse additive models is spirit of \citet
{MGB09}. More precisely,
if ${\mathcal H}_j:=\{h\in W^{\alpha,2}[0,1]\dvtx\int_0^1 h(x)\,dx=0\}$,
then, for uniform distribution $\Pi$, the spaces ${\mathcal H}_j$ are
orthogonal in $L_2(\Pi)$ (recall that ${\mathcal H}_j$ is viewed as
a space of functions depending on the $j$th coordinate). Assume, for
simplicity, that $\ell$ is the quadratic loss and that the regression
function $f_{\ast}$ can be represented as $f_{\ast}=\sum_{j\in
J}f_{\ast
,j}$, where $J$ is a subset of $\{1,\ldots, N\}$ of cardinality
$d$ and $\|f_{\ast,j}\|_{{\mathcal H}_j}\leq1$. Then it easily follows
from the bound of Theorem~\ref{co:homobd-1} that with probability at least
$1-3N^{-A/2}$
\[
{\mathcal E}(f)=\|f-f_{\ast}\|_{L_2(\Pi)}^2
\leq C\tau^2 d \biggl(n^{-2\alpha/(2\alpha+1)}\vee\frac{A\log
N}{n} \biggr).
\]
Note that, up to a constant, this essentially coincides with the
minimax lower bound in this type of problems obtained recently by
\citet{RaskuttiWainwrightYu2009}.
Of course, if the design
distribution is not necessarily uniform, an adaptive choice
of regularization parameters might be needed even in such
simple examples and the approach described above leads to
minimax optimal rates.
\end{Remarks*}

\section{Preliminary bounds}
\label{sec:sam2pop}

In this section, the case of a single RKHS $\mathcal{H}_K$ associated with
a kernel $K$ is considered. We assume that $K(x,x)\leq1, x\in S$.
This implies
that, for all $h\in\Hcal_K$,
$
\|h\|_{L_2(\Pi)}\le\|h\|_{L_{\infty}}\leq\|h\|_{\Hcal_K}.
$

\subsection{\texorpdfstring{Comparison of $\|\cdot\|_{L_2(\Pi_n)}$ and $\|\cdot\|_{L_2(\Pi)}$}
{Comparison of ||.|| L 2(Pi n) and ||.|| L 2(Pi)}}

First, we study the relationship between the empirical and
the population $L_2$ norms
for functions in $\Hcal_K$.
\begin{theorem}
\label{th:empnorm}
Assume that $A\geq1$ and $\log N\geq2\log\log n$.
Then there exists a numerical constant $C>0$ such that
with probability at least $1-N^{-A}$ for all $h\in\Hcal_K$
%
%
\begin{eqnarray}
\label{eq:empl2-1}
\|h\|_{L_2(\Pi)}&\le& C \bigl(\|h\|_{L_2(\Pi_n)}+
\bar{\epsilon}\|h\|_{\Hcal_K} \bigr);\\
\label{eq:empl2-2}
\|h\|_{L_2(\Pi_n)}&\le& C \bigl(\|h\|_{L_2(\Pi)}+
\bar{\epsilon}\|h\|_{\Hcal_K} \bigr),
\end{eqnarray}
where
%
%
\begin{eqnarray}\quad
\bar{\epsilon}&=&\bar{\epsilon}(K)\nonumber\\[-8pt]\\[-8pt]
:\!&=&
\inf\Biggl\{\epsilon\ge\sqrt{A\log N\over n}\dvtx
{\bE\mathop{\sup_{\|h\|_{\Hcal_K}=1}}_{\|h\|_{L_2(\Pi)}\le\delta}}
|R_n (h)|\le\epsilon\delta+ \epsilon^2, \forall\delta\in(0,1]\Biggr\}.\nonumber
\end{eqnarray}
\end{theorem}
\begin{pf}
Observe that the inequalities
hold trivially when $h=0$.
We shall therefore consider only the case when $h\neq0$.
By symmetrization inequality,
%
%
\begin{equation}\qquad
{\bE\mathop{\sup_{\|h\|_{\Hcal_K}=1}}_{2^{-j}<
\|h\|_{L_2(\Pi)}\le2^{-j+1}}} |(\Pi_n-\Pi)h^2 |
\le{2\bE\mathop{\sup_{\|h\|_{\Hcal_K}=1}}_{2^{-j}<\|h\|_{L_2(\Pi
)}\le2^{-j+1}}}
|R_n (h^2)|
\end{equation}
and, by contraction inequality, we further have
%
%
\begin{equation}\qquad
{\bE\mathop{\sup_{\|h\|_{\Hcal_K}=1}}_{2^{-j}<\|h\|_{L_2(\Pi)}
\le2^{-j+1}}} |(\Pi_n-\Pi)h^2 |\le
{8\bE\mathop{\sup_{\|h\|_{\Hcal_K}=1}}_{2^{-j}<\|h\|_{L_2(\Pi)}
\le2^{-j+1}}} |R_n (h) |.
\end{equation}

The definition of $\bar{\epsilon}$ implies that
%
%
\begin{eqnarray}
&&{\bE\mathop{\sup_{\|h\|_{\Hcal_K}=1}}_{2^{-j}<\|h\|_{L_2(\Pi)}
\le2^{-j+1}}} |(\Pi_n-\Pi)h^2 |\nonumber\\[-8pt]\\[-8pt]
&&\qquad\le
{8 \bE\mathop{\sup_{\|h\|_{\Hcal_K}=1}}_{\|h\|_{L_2(\Pi)}\le
2^{-j+1}}} |R_n(h) |
\le8 (\bar{\epsilon}2^{-j+1}+\bar{\epsilon}^2 ).\nonumber
\end{eqnarray}

An application of Talagrand's concentration inequality yields
\begin{eqnarray*}
&&
{\mathop{\sup_{\|h\|_{\Hcal_K}=1}}_{2^{-j}<\|h\|_{L_2(\Pi)}
\le2^{-j+1}}} |(\Pi_n-\Pi)h^2 |\\
&&\qquad\le
2 \Biggl({\bE\mathop{\sup_{\|h\|_{\Hcal_K}=1}}_{2^{-j}<\|h\|
_{L_2(\Pi)}
\le2^{-j+1}}} |(\Pi_n-\Pi)h^2 |\\
&&\qquad\quad\hspace*{9.2pt}{}+2^{-j+1}\sqrt{t+2\log j\over n}+{t+2\log j \over
n}
\Biggr)\\
&&\qquad\leq 32 \Biggl(\bar{\epsilon}2^{-j}+\bar{\epsilon}^2+
2^{-j}\sqrt{t+2\log j\over n}+{t+2\log j \over n} \Biggr)
\end{eqnarray*}
with probability at least $1-\exp(-t-2\log j)$ for any natural number $j$.
Now, by the union bound, for all $j$ such that $2\log j\le t$,
%
%
\begin{eqnarray}\label{Pi_n-Pi}
&&{\mathop{\sup_{\|h\|_{\Hcal_K}=1}}_{2^{-j}<\|h\|_{L_2(\Pi)}\le
2^{-j+1}}} |(\Pi_n-\Pi)h^2 |\nonumber\\[-8pt]\\[-8pt]
&&\qquad
\le32 \Biggl(\bar{\epsilon}2^{-j}+\bar{\epsilon}^2+
2^{-j}\sqrt{t+2\log j\over n}+{t+2\log j \over n} \Biggr)\nonumber
\end{eqnarray}
with probability at least
%
%
\begin{eqnarray}
1-\sum_{j\dvtx2\log j\le t} \exp(-t -2\log j)&=&
1-\exp(-t)\sum_{j\dvtx2\log j\le t}j^{-2}\nonumber\\[-8pt]\\[-8pt]
&\ge&1-2\exp(-t).\nonumber
\end{eqnarray}
Recall that $\bar{\epsilon}\geq(A\log N/n)^{1/2}$
and $\|h\|_{L_2(\Pi)}\le\|h\|_{\Hcal_K}$.
Taking $t=A\log N+\log4$, we easily get that,
for all $h\in\mathcal{H}_K$
such that $\|h\|_{\mathcal{H}_K}=1$ and $\|h\|_{L_2(\Pi)}\geq
\exp\{-N^{A/2}\}$,
%
%
\begin{equation}
|(\Pi_n-\Pi)h^2 |\le
C \bigl({\bar{\epsilon}\|h\|_{L_2(\Pi)}} +
\bar{\epsilon}^2 \bigr)
\end{equation}
with probability at least $1-0.5 N^{-A}$ and with a numerical
constant $C>0$.
In other words, with the same
probability, for all $h\in\mathcal{H}_K$ such that
$\frac{\|h\|_{L_2(\Pi)}}{\|h\|_{\mathcal{H}_K}}\geq
\exp\{-N^{A/2}\}$,
%
%
\begin{equation}
|(\Pi_n-\Pi)h^2 |
\le C \bigl(\bar{\epsilon}\|h\|_{L_2(\Pi)}\|h\|_{\Hcal_K}
+\bar{\epsilon}^2\|h\|_{\Hcal_K}^2 \bigr).
\end{equation}
Therefore, for all $h\in\Hcal_K$ such that
%
%
\begin{equation}
{\|h\|_{L_2(\Pi)}\over\|h\|_{\Hcal_K}}>\exp(-N^{A/2})
\end{equation}
we have
\begin{eqnarray*}
\|h\|_{L_2(\Pi)}^2&=&
\Pi h^2 \le \|h\|_{L_2(\Pi_n)}^2+
C \bigl(\bar{\epsilon}\|h\|_{L_2(\Pi)}\|h\|_{\Hcal_K}
+\bar{\epsilon}^2\|h\|_{\Hcal_K}^2 \bigr), \\
\|h\|_{L_2(\Pi_n)}^2&=&\Pi_n h^2\le
\|h\|_{L_2(\Pi)}^2+C \bigl(\bar{\epsilon}\|h\|_{L_2(\Pi)}
\|h\|_{\Hcal_K} +\bar{\epsilon}^2\|h\|_{\Hcal_K}^2 \bigr).
\end{eqnarray*}
It can be now deduced that, for a proper value of numerical constant $C$,
%
%
\begin{eqnarray}
\|h\|_{L_2(\Pi)} &\le& C \bigl(\|h\|_{L_2(\Pi_n)}+
\bar{\epsilon}\|h\|_{\Hcal_K} \bigr)\quad
\mbox{and}\nonumber\\[-8pt]\\[-8pt]
\|h\|_{L_2(\Pi_n)} &\le& C \bigl(\|h\|_{L_2(\Pi)}+
\bar{\epsilon}\|h\|_{\Hcal_K} \bigr).\nonumber
\end{eqnarray}

It remains to consider the case when
%
%
\begin{equation}
{\|h\|_{L_2(\Pi)}\over\|h\|_{\Hcal_K}}\le\exp(-N^{A/2}).
\end{equation}
Following a similar argument as before, with probability at least $1-0.5
N^{-A}$,
\begin{eqnarray*}
&&{\mathop{\sup_{\|h\|_{\Hcal_K}=1}}_{\|h\|_{L_2(\Pi)}
\le\exp(-N^{A/2})}}
|(\Pi_n-\Pi)h^2 |\\
&&\qquad\le
16 \Biggl(\bar{\epsilon}\exp(-N^{A/2})+\bar{\epsilon}^2
+\exp(-N^{A/2})\sqrt{A\log N\over n}+{A\log N \over
n} \Biggr).
\end{eqnarray*}
Under the conditions $A\geq1, \log N\geq2 \log\log n$,
%
%
\begin{equation}
\bar{\epsilon}\ge\biggl({A\log N\over n} \biggr)^{1/2}\geq\exp(-N^{A/2}).
\end{equation}
Then
%
%
\begin{equation}
{\mathop{\sup_{\|h\|_{\Hcal_K}=1}}_{\|h\|_{L_2(\Pi)}
\le\exp(-N^{A/2})} }|(\Pi_n-\Pi)h^2 |
\le C\bar{\epsilon}^2
\end{equation}
with probability at least $1-0.5 N^{-A}$,
which also implies (\ref{eq:empl2-1}) and (\ref{eq:empl2-2}),
and the result follows.
\end{pf}

Theorem~\ref{th:empnorm} shows that the two norms
$\|h\|_{L_2(\Pi_n)}$ and $\|h\|_{L_2(\Pi)}$
are of the same order up to an error term
$\bar{\epsilon}\|h\|_{\Hcal_K}$.

\subsection{\texorpdfstring{Comparison of $\hat{\epsilon}(K)$, $\bar{\epsilon}(K)$, $\breve\epsilon(K)$ and $\check{\epsilon}(K)$}
{Comparison of epsilon(K), epsilon(K), epsilon(K) and epsilon(K)}}

Recall the definitions
\[
\breve\gamma_n (\delta):=
\Biggl(n^{-1}\sum_{k=1}^\infty
(\lambda_k\wedge\delta^2) \Biggr)^{1/2},\qquad \delta\in(0,1],
\]
where $\{\lambda_k\}$ are the eigenvalues of the integral operator $T_K$
from $L_2(\Pi)$ into $L_2(\Pi)$ with kernel $K$, and, for some $A\geq1$,
\[
\breve{\epsilon}(K):=
\inf\Biggl\{\epsilon\ge\sqrt{A\log N\over n}\dvtx
\breve\gamma_n(\delta)\le\epsilon\delta+ \epsilon^2,
\forall\delta\in(0,1] \Biggr\}.
\]
It follows from Lemma 42 of \citet{Mendelson2002}
[with an additional application of Cauchy--Schwarz inequality
for the upper bound and Hoffmann--J\o rgensen inequality
for the lower bound; see also \citet{Koltchinskii2008}]
that, for some
numerical constants $C_1,C_2>0$,
%
%
\begin{eqnarray}
\label{thirtythree}
C_1
\Biggl(n^{-1}\sum_{k=1}^n (\lambda_k \wedge\delta^2) \Biggr)^{1/2}-n^{-1}
&\leq&
{\bE
\mathop{\sup_{\|h\|_{\Hcal_K}=1}}_{\|h\|_{L_2(\Pi)}
\le\delta}} |R_n (h) |\nonumber\\[-8pt]\\[-8pt]
&\leq&
C_2
\Biggl(n^{-1}\sum_{k=1}^n (\lambda_k \wedge\delta^2)
\Biggr)^{1/2}.\nonumber
\end{eqnarray}

This fact and the definitions of $\breve\epsilon(K), \bar\epsilon(K)$
easily imply the following result.
\begin{proposition}
\label{pr:bareps}
Under the condition $K(x,x)\leq1, x\in S$,
there exist numerical constants $C_1, C_2>0$ such that
%
%
\begin{equation}
C_1\breve\epsilon(K)
\le\bar{\epsilon}(K)\le
C_2\breve\epsilon(K).
\end{equation}
\end{proposition}

If $K$ is the kernel of the projection operator onto a finite-dimensional
subspace $\mathcal{H}_K$ of $L_2(\Pi)$, it is easy to check that
$
\breve\epsilon(K)\asymp\sqrt{\frac{\operatorname{dim}(\mathcal{H}_K)}{n}}
$
(recall the notation $a\asymp b$, which means
that there exists a numerical constant $c>0$ such that
\mbox{$c^{-1}\leq{a}/{b}\leq c$}). If the eigenvalues $\lambda_k$ decay at a
polynomial rate, that is,
$\lambda_k \asymp k^{-2\beta}$
for some $\beta>1/2$, then
$\breve\epsilon(K)\asymp n^{-\beta/(2\beta+1)}$.

Recall the notation
%
%
\begin{equation}
\hat{\epsilon}(K):=
\inf\Biggl\{\epsilon\ge\sqrt{A\log N\over n}\dvtx
\Biggl({1\over n}\sum_{k=1}^n (\hat\lambda_k\wedge
\delta^2 ) \Biggr)^{1/2}\le\epsilon\delta+ \epsilon^2,
\forall\delta\in(0,1] \Biggr\},\hspace*{-35pt}
\end{equation}
where $\{\hat\lambda_k\}$ denote the eigenvalues of the Gram matrix
$
\hat K:= (K(X_i,X_j) )_{i,j=1,\ldots, n}.
$
It follows again from the results of \citet{Mendelson2002}
[namely, one can follow the proof of Lemma 42 in the case
when the RKHS $\mathcal{H}_K$ is restricted to the sample $X_1,\ldots, X_n$
and the expectations are conditional on the sample; then one uses
Cauchy--Schwarz and Hoffmann--J\o rgensen inequalities as in the proof of
(\ref{thirtythree})] that for some numerical
constants $C_1,C_2>0$
%
%
\begin{eqnarray}
C_1
\Biggl(n^{-1}\sum_{k=1}^n (\hat\lambda_k \wedge\delta^2)
\Biggr)^{1/2}-n^{-1}
&\leq&
{\bE_{\varepsilon}
\mathop{\sup_{\|h\|_{\Hcal_K}=1}}_{\|h\|_{L_2(\Pi_n)}
\le\delta}} |R_n (h) |\nonumber\\[-8pt]\\[-8pt]
&\leq&
C_2
\Biggl(n^{-1}\sum_{k=1}^n (\hat\lambda_k \wedge\delta^2) \Biggr)^{1/2},\nonumber
\end{eqnarray}
where $\bE_{\varepsilon}$ indicates that the expectation
is taken over the Rademacher random variables only (conditionally
on $X_1,\ldots, X_n$).
Therefore, if we denote~by
%
%
\begin{equation}
\tilde{\epsilon}(K):=
\inf\Biggl\{\epsilon\ge\sqrt{A\log N\over n}\dvtx
\bE_{\varepsilon}\mathop{\sup_{\|h\|_{\Hcal_K}=1}}_{
\|h\|_{L_2(\Pi_n)}\le\delta}|R_n(h)|
\le\epsilon\delta+ \epsilon^2, \forall\delta\in(0,1] \Biggr\}\hspace*{-35pt}
\end{equation}
the empirical version of $\bar{\epsilon}(K)$,
then $\hat\epsilon(K)\asymp\tilde\epsilon(K)$.
We will now show that $\tilde\epsilon(K) \asymp\bar\epsilon(K)$
with a high probability.
\begin{theorem}
\label{th:empeps}
Suppose that $A\geq1$ and $\log N\geq2\log\log n$.
There exist numerical constants $C_1, C_2>0$ such that
%
%
\begin{equation}
C_1\bar{\epsilon}(K)\le\tilde{\epsilon}(K)\le C_2\bar{\epsilon}(K),
\end{equation}
with probability at least $1-N^{-A}$.
\end{theorem}
\begin{pf}
Let $t:=A\log N+\log14$.
It follows from
Talagrand concentration inequality that
\begin{eqnarray*}
&&{\bE\mathop{\sup_{\|h\|_{\Hcal_K}=1}}_{2^{-j}<\|h\|_{L_2(\Pi)}\le
2^{-j+1}}}
|R_n(h) |\\
&&\qquad\le 2 \Biggl({\mathop{\sup_{\|h\|_{\Hcal_K}=1}}_{
2^{-j}<\|h\|_{L_2(\Pi)}\le2^{-j+1}}}
|R_n(h) |+2^{-j+1}\sqrt{t+2\log j\over n}+{t+2\log j\over
n} \Biggr)
\end{eqnarray*}
with probability at least $1-\exp(-t-2\log j)$.
On the other hand, as derived in the proof of
Theorem~\ref{th:empnorm} [see (\ref{Pi_n-Pi})]
%
%
\begin{eqnarray}
&&
{\mathop{\sup_{\|h\|_{\Hcal_K}=1}}_{2^{-j}<\|h\|_{L_2(\Pi)}\le2^{-j+1}}}
|(\Pi_n-\Pi)h^2 |\nonumber\\[-8pt]\\[-8pt]
&&\qquad\le32 \Biggl(\bar{\epsilon}2^{-j}+\bar
{\epsilon}^2+
2^{-j}\sqrt{t+2\log j\over n}+{t+2\log j \over n} \Biggr)\nonumber
\end{eqnarray}
with probability at least $1-\exp(-t-2\log j)$.
We will use these bounds only for $j$ such that $2\log j\leq t$.
In this case, the second bound implies that, for some numerical
constant $c>0$ and all $h$ satisfying the conditions
$
\|h\|_{\Hcal_K}=1, 2^{-j}<\|h\|_{L_2(\Pi)}\le2^{-j+1},
$
we have
$
\|h\|_{L_2(\Pi_n)}\leq c(2^{-j}+\bar\epsilon)
$
(again, see the proof of
Theorem~\ref{th:empnorm}).
Combining these bounds, we get that with probability
at least $1-2\exp(-t-2\log j)$,
\begin{eqnarray*}
&&
{\bE\mathop{\sup_{\|h\|_{\Hcal_K}=1}}_{2^{-j}<\|h\|_{L_2(\Pi)}\le2^{-j+1}}}
|R_n(h) |\\
&&\qquad\le2 \Biggl({\mathop{\sup_{\|h\|_{\Hcal_K}=1}}_{
\|h\|_{L_2(\Pi_n)}\le c\delta_j}}
|R_n(h) |+2^{-j+1}\sqrt{t+2\log j\over n}+{t+2\log j\over
n} \Biggr),
\end{eqnarray*}
where $\delta_j=\bar{\epsilon}+2^{-j}$.

Applying now Talagrand concentration inequality
to the Rademacher process
conditionally on the observed data $X_1,\ldots, X_n$
yields
\begin{eqnarray*}
&&{\mathop{\sup_{\|h\|_{\Hcal_K}=1}}_{
\|h\|_{L_2(\Pi_n)}\le c\delta_j}}
|R_n(h) |
\le
2 \Biggl({\bE_\varepsilon\mathop{\sup_{\|h\|_{\Hcal_K}=1}}_{
\|h\|_{L_2(\Pi_n)}\le c\delta_j}}
|R_n(h) |\\
&&\qquad\quad\hspace*{79.3pt}{}
+C\delta_j\sqrt{t+2\log j\over n}+{t+2\log j\over n} \Biggr),
\end{eqnarray*}
with conditional probability at least $1-\exp(-t-2\log j)$.
From this and from the previous bound it is not hard to deduce that,
for some numerical constants $C, C^{\prime}$ and for all $j$ such
that $2\log j\leq t$,
\begin{eqnarray*}
&&{\bE\mathop{\sup_{\|h\|_{\Hcal_K}=1}}_{2^{-j}<\|h\|_{L_2(\Pi)}\le
2^{-j+1}}}
|R_n(h) |\\
&&\qquad\le C' \Biggl({\bE_\varepsilon\mathop{\sup_{\|h\|_{\Hcal
_K}=1}}_{\|h\|
_{L_2(\Pi_n)}\le c\delta_j}}
|R_n(h) |+\delta_j\sqrt{t+2\log j\over n}+
{t+2\log j\over n} \Biggr)\\
&&\qquad\le C(\tilde{\epsilon}\delta_j+\tilde{\epsilon}^2)
\le C(\tilde\epsilon2^{-j} + \tilde\epsilon\bar\epsilon+
\tilde\epsilon^2)
\end{eqnarray*}
with probability at least $1-3\exp(-t-2\log j)$.
In obtaining the second inequality, we used the definition of
$\tilde{\epsilon}$ and the fact that, for $t=A\log N+\log14, 2\log
j\leq t$,
$c_1\tilde\epsilon\geq(t+2\log j/n)^{1/2}$, where $c_1$ is a numerical
constant. Now, by the union bound,
the above inequality holds with probability at least
%
%
\begin{equation}
1-3\sum_{j\dvtx2\log j\le t}\exp(-t-2\log j)\ge1-6\exp(-t)
\end{equation}
for all $j$ such that $2\log j\le t$ simultaneously.
Similarly, it can be shown that
\[
{\bE\mathop{\sup_{\|h\|_{\Hcal_K}=1}}_{\|h\|_{L_2(\Pi)}\le\exp(-N^{A/2})}}
|R_n(h) |
\le
C \bigl(\tilde{\epsilon}\exp(-N^{A/2})+
\tilde{\epsilon}\bar{\epsilon}+
\tilde{\epsilon}^2 \bigr)
\]
with probability at least $1-\exp(-t)$.

For $t=A\log N+\log14$, we get
%
%
\begin{equation}
{\bE\mathop{\sup_{\|h\|_{\Hcal_K}=1}}_{\|h\|_{L_2(\Pi)}\le\delta}}
|R_n(h) |\le
C (\tilde{\epsilon}\delta+ \tilde{\epsilon}\bar{\epsilon} +
\tilde{\epsilon}^2 ),
\end{equation}
for all $0<\delta\le1$,
with probability at least
$1-7 \exp(-t)=
1-N^{-A}/2.
$
Now by the definition of $\bar{\epsilon}$, we obtain
%
%
\begin{equation}
\bar{\epsilon}\le C\max\{\tilde{\epsilon},(\tilde{\epsilon}\bar
{\epsilon}
+\tilde{\epsilon}^2)^{1/2}\},
\end{equation}
which implies that $\bar{\epsilon}\le C\tilde{\epsilon}$
with probability at least $1-N^{-A}/2$.

Similarly one can show that
%
%
\begin{equation}
{\bE_\epsilon\mathop{\sup_{\|h\|_{\Hcal_K}=1}}_{\|h\|_{L_2(\Pi
)}\le
\delta}}
|R_n(h) |\le
C (\bar{\epsilon}\delta+ \tilde{\epsilon}\bar{\epsilon} +
\bar{\epsilon}^2 ),
\end{equation}
for all $0<\delta\le1$, with probability at least
$1-N^{-A}/2$,
which implies that $\tilde{\epsilon}\le C\bar{\epsilon}$ with
probability at least $1-N^{-A}/2$. The proof can then be completed by
the union bound.
\end{pf}

Define
%
%
\begin{eqnarray}
\check\epsilon &:=& \check{\epsilon}(K)\nonumber\\[-8pt]\\[-8pt]
&:=&
\inf\Biggl\{\epsilon\ge\sqrt{A\log N\over n}\dvtx
\mathop{\sup_{\|h\|_{{\Hcal}_K}=1}}_{
\|h\|_{L_2(\Pi)}
\le\delta}|R_n(h)|\le
\epsilon\delta+ \epsilon^2, \forall\delta\in(0,1] \Biggr\}.\nonumber
\end{eqnarray}

The next statement can be proved similarly to Theorem~\ref{th:empeps}.
\begin{theorem}
\label{th:empeps-1}
There exist numerical constants $C_1, C_2>0$ such that
%
%
\begin{equation}
C_1\bar{\epsilon}(K)\le\check{\epsilon}(K)\le C_2\bar{\epsilon}(K)
\end{equation}
with probability at least $1-N^{-A}$.
\end{theorem}

Suppose now that $\{K_1,\ldots, K_N\}$ is a dictionary of kernels.
Recall that $\bar\epsilon_j=\bar\epsilon(K_j)$,
$\hat\epsilon_j=\hat\epsilon(K_j)$ and
$\check\epsilon_j=\check\epsilon(K_j)$.

It follows from Theorems~\ref{th:empnorm},~\ref{th:empeps},
\ref{th:empeps-1} and the union bound that
with probability at least $1-3N^{-A+1}$ for all $j=1,\ldots, N$
%
%
\begin{eqnarray}
\label{conclusion_A}
\|h\|_{L_2(\Pi)}&\le& C \bigl(\|h\|_{L_2(\Pi_n)}+
\bar{\epsilon}_j\|h\|_{\Hcal_K} \bigr),\nonumber\\[-8pt]\\[-8pt]
\|h\|_{L_2(\Pi_n)}&\le& C \bigl(\|h\|_{L_2(\Pi)}+
\bar{\epsilon}_j\|h\|_{\Hcal_K} \bigr),\qquad h\in\mathcal{H}_j,\nonumber
%
%
\\
\label{conclusion_B}
C_1\bar{\epsilon}_j&\le&\hat{\epsilon}_j\le C_2\bar{\epsilon}_j
\quad\mbox{and}\quad
C_1\bar{\epsilon}_j\le\check{\epsilon}_j\le C_2\bar{\epsilon}_j.
\end{eqnarray}
Note also that
\[
3N^{-A+1}=\exp\{-(A-1)\log N+\log3\}\leq\exp\{-(A/2)\log N\}=N^{-A/2},
\]
provided that $A\geq4$ and $N\geq3$. Thus, under these additional
constraints, (\ref{conclusion_A}) and (\ref{conclusion_B}) hold for all
$j=1,\ldots, N$ with probability at least $1-N^{-A/2}$.

\section{Proofs of the oracle inequalities}
\label{sec:mainres}

For an arbitrary set $J\subseteq\{1,\ldots, N\}$ and $b\in(0,+\infty)$,
denote
%
%
\begin{equation}\quad
\mathcal{K}_J^{(b)}:= \biggl\{(f_1,\ldots, f_N)\in\Hcal^{(N)}\dvtx
\sum_{j\notin J}
\bar{\epsilon}_j\|f_j\|_{L_2(\Pi)}
\le b\sum_{j\in J}\bar{\epsilon}_j \|f_j\|_{L_2(\Pi)} \biggr\}
\end{equation}
and let
%
%
\begin{eqnarray}
\beta_b(J)&=&\inf\biggl\{{\beta}\geq0\dvtx\sum_{j\in J}
\bar{\epsilon}_j\|f_j\|_{L_2(\Pi)}\le
\beta\|f_1+\cdots+f_N\|_{L_2(\Pi)},\nonumber\\[-8pt]\\[-8pt]
&&\hspace*{153pt} (f_1,\ldots, f_N)\in\mathcal
{K}_J^{(b)} \biggr\}.\nonumber
\end{eqnarray}

It is easy to see that, for all nonempty sets $J$,
$
\beta_b(J)\geq\max_{j\in J}\bar\epsilon_j \geq\sqrt{\frac{A\log N}{n}}.
$

Theorems~\ref{co:homobd} and~\ref{co:homobd-1} will be easily
deduced from the following technical result.
\begin{theorem}
\label{th:main}
There exist numerical constants $C_1,C_2,B>0$ and $b>0$ such that, for all
$\tau\ge B L_{\ast}$ in the definition of
$\epsilon_j=\tau\hat\epsilon_j, j=1,\ldots, N$
and for all oracles $(f_1,\ldots, f_N) \in D$,
%
%
\begin{eqnarray}
\label{main_oracle}
&&\Ecal(\ell\circ\hat{f})+C_1 \Biggl(\sum_{j=1}^N
\tau\bar{\epsilon}_j\|\hat{f}_j-f_j\|_{L_2(\Pi)}
+\sum_{j=1}^N\tau^2\bar{\epsilon}_j^2\|\hat{f}_j\|_{\Hcal
_j} \Biggr)\\
&&\qquad\le
2\Ecal(\ell\circ f)+C_2\tau^2 \biggl(\sum_{j\in J_f}\bar{\epsilon
}_j^2\|
f_j\|_{\Hcal_j}+\frac{\beta_b^2(J_f)}{m_{\ast}} \biggr)
\end{eqnarray}
with probability at least $1-3 N^{-A/2}$.
Here, $A\geq4$ is a constant involved in the definitions
of $\bar\epsilon_j, \hat\epsilon_j, j=1,\ldots, N$.
\end{theorem}
\begin{pf}
Recall that
\begin{eqnarray*}
(\hat{f}_1,\ldots, \hat{f}_N )&:=&
\argmin_{(f_1,\ldots, f_N)\in D}
\Biggl[P_n\bigl(\ell\circ(f_1+\cdots+f_N)\bigr)\\
&&\hspace*{54.4pt}{}+\sum_{j=1}^N
\bigl(\tau\hat{\epsilon}_j\|f_j\|_{L_2(\Pi_n)} +
\tau^2\hat{\epsilon}_j^2\|f_j\|_{\Hcal_j} \bigr) \Biggr],
\end{eqnarray*}
and that we write
$
f:=f_1+\cdots+f_N, \hat f:=\hat f_1+\cdots+\hat f_N.
$
Hence, for all $(f_1,\ldots, f_N)\in D$,
\begin{eqnarray*}
&& P_n(\ell\circ\hat f)
+\sum_{j=1}^N \bigl(\tau\hat{\epsilon}_j\|\hat{f}_j\|_{L_2(\Pi_n)}
+\tau^2\hat{\epsilon}_j^2\|\hat{f}_j\|_{\Hcal_j} \bigr)\\
&&\qquad\le P_n(\ell\circ f)
+\sum_{j=1}^N \bigl(\tau\hat{\epsilon}_j\|f_j\|_{L_2(\Pi_n)}
+\tau^2\hat{\epsilon}_j^2\|f_j\|_{\Hcal_j} \bigr).
\end{eqnarray*}
By a simple algebra,
\begin{eqnarray*}
&&\Ecal(\ell\circ\hat{f})+\sum_{j=1}^N
\bigl(\tau\hat{\epsilon}_j\|\hat{f}_j\|_{L_2(\Pi_n)}
+\tau^2\hat{\epsilon}_j^2\|\hat{f}_j\|_{\Hcal_j} \bigr)\\
&&\qquad\le\Ecal(\ell\circ f)+\sum_{j=1}^N
\bigl(\tau\hat{\epsilon}_j\|f_j\|_{L_2(\Pi_n)}+
\tau^2\hat{\epsilon}_j^2\|f_j\|_{\Hcal_j} \bigr)\\
&&\qquad\quad{}
+ |(P_n-P) (\ell\circ\hat{f}-\ell\circ f ) |
\end{eqnarray*}
and, by the triangle inequality,
\begin{eqnarray*}
&&\Ecal(\ell\circ\hat{f})+
\sum_{j\notin J_f}\tau\hat{\epsilon}_j\|\hat{f}_j\|_{L_2(\Pi_n)}
+\sum_{j=1}^N\tau^2\hat{\epsilon}_j^2\|\hat{f}_j\|_{\Hcal_j}\\
&&\qquad\le\Ecal(\ell\circ f)+\sum_{j\in J_f}\tau\hat{\epsilon}_j
\|\hat{f}_j-f_j\|_{L_2(\Pi_n)}\\
&&\qquad\quad{} +
\sum_{j\in J_f}\tau^2\hat{\epsilon}_j^2\|f_j\|_{\Hcal_j}
+ |(P_n-P) (\ell\circ\hat{f}-\ell\circ f ) |.
\end{eqnarray*}

We now take advantage of (\ref{conclusion_A}) and (\ref{conclusion_B})
to replace $\hat{\epsilon}_j$'s by $\bar{\epsilon}_j$'s
and $\|\cdot\|_{L_2(\Pi_n)}$ by $\|\cdot\|_{L_2(\Pi)}$.
Specifically, there exists a numerical constant $C>1$ and an event $E$
of probability
at least $1-N^{-A/2}$ such that
%
%
\begin{equation}
\label{condition_01}
{1\over C}\le\min\biggl\{{\hat{\epsilon}_j\over\bar{\epsilon}_j}\dvtx
j=1,\ldots, N \biggr\}\le\max\biggl\{{\hat{\epsilon}_j\over\bar
{\epsilon
}_j}\dvtx j=1,\ldots, N \biggr\}\le C
\end{equation}
and, for all $j=1,\ldots,N$,
%
%
\begin{equation}\quad
\label{condition_02}
{1\over C}\|\hat{f}_j\|_{L_2(\Pi)}-
\bar{\epsilon}_j\|\hat{f}_j\|_{\Hcal_j}
\le\|\hat{f}_j\|_{L_2(\Pi_n)}\le
C \bigl(\|\hat{f}_j\|_{L_2(\Pi)}
+\bar{\epsilon}_j\|\hat{f}_j\|_{\Hcal_j} \bigr).
\end{equation}
Taking $\tau\geq C/(C-1)$, we have that, on the event $E$,
\begin{eqnarray*}
\hspace*{-5pt}&&\Ecal(\ell\circ\hat{f})+\sum_{j\notin J_f}
\tau\hat{\epsilon}_j\|\hat{f}_j\|_{L_2(\Pi_n)}
+\sum_{j=1}^N\tau^2\hat{\epsilon}_j^2\|\hat{f}_j\|_{\Hcal_j}\\
\hspace*{-5pt}&&\qquad\ge \Ecal(\ell\circ\hat{f})
+{1\over C^2}
\Biggl(\sum_{j\notin J_f}\tau\bar{\epsilon}_j\|\hat{f}_j\|
_{L_2(\Pi_n)}
+\sum_{j=1}^N\tau^2\bar{\epsilon}_j^2\|\hat{f}_j\|_{\Hcal
_j} \Biggr)\\
\hspace*{-5pt}&&\qquad\ge
\Ecal(\ell\circ\hat{f})+
{1\over C^2}
\Biggl(\sum_{j\notin J_f}\tau\bar{\epsilon}_j
\biggl({1\over C}\|\hat{f}_j\|_{L_2(\Pi)}
-\bar{\epsilon}_j\|\hat{f}_j\|_{\Hcal_j} \biggr)
+\sum_{j=1}^N\tau^2\bar{\epsilon}_j^2\|\hat{f}_j\|_{\Hcal
_j} \Biggr)\\
\hspace*{-5pt}&&\qquad\ge \Ecal(\ell\circ\hat{f})+
{1\over C^3}
\Biggl(\sum_{j\notin J_f}\tau\bar{\epsilon}_j
\|\hat{f}_j\|_{L_2(\Pi)}+
\sum_{j=1}^N\tau^2\bar{\epsilon}_j^2\|\hat{f}_j\|_{\Hcal_j} \Biggr).
\end{eqnarray*}
Similarly,
\begin{eqnarray*}
&&\Ecal(\ell\circ f)+\sum_{j\in J_f}
\bigl(\tau\hat{\epsilon}_j\|f_j-\hat{f}_j\|_{L_2(\Pi_n)}
+\tau^2\hat{\epsilon}_j^2\|f_j\|_{\Hcal_j} \bigr)\\
&&\qquad\le
\Ecal(\ell\circ f)+C^2\sum_{j\in J_f}
\bigl(\tau\bar{\epsilon}_j\|f_j-\hat{f}_j\|_{L_2(\Pi_n)}
+\tau^2\bar{\epsilon}_j^2\|f_j\|_{\Hcal_j} \bigr)\\
&&\qquad\le\Ecal(\ell\circ f)+
C^3\sum_{j\in J_f}\tau\bar{\epsilon}_j
\bigl(\|{f}_j-\hat{f}_j\|_{L_2(\Pi)}
+\bar{\epsilon}_j\|f_j-\hat{f}_j\|_{\Hcal_j} \bigr)\\
&&\qquad\quad{}+C^2\sum_{j\in J_f}\tau^2\bar{\epsilon}_j^2\|f_j\|_{\Hcal_j}\\
&&\qquad\le\Ecal(\ell\circ f)+C^3\sum_{j\in J_f}\tau\bar{\epsilon}_j
\bigl(\|{f}_j-\hat{f}_j\|_{L_2(\Pi)}
+\bar{\epsilon}_j\|f_j\|_{\Hcal_j}
+\bar{\epsilon}_j\|\hat{f}_j\|_{\Hcal_j} \bigr)\\
&&\qquad\quad{}+C^2\sum_{j\in J_f}\tau^2\bar{\epsilon}_j^2\|f_j\|_{\Hcal_j}\\
&&\qquad\le\Ecal(\ell\circ f)+
2C^3\sum_{j\in J_f} \bigl(\tau\bar{\epsilon}_j
\|{f}_j-\hat{f}_j\|_{L_2(\Pi)}
+\tau^2\bar{\epsilon}_j^2\|f_j\|_{\Hcal_j} \bigr)\\
&&\qquad\quad{}+
C^3\sum_{j\in J_f}\tau\bar{\epsilon}_j^2\|\hat{f}_j\|_{\Hcal_j}.
\end{eqnarray*}
Therefore, by taking $\tau$ large enough, namely
$\tau\geq\frac{C}{C-1}\vee(2C^6)$,
we can find numerical constants $0<C_1<1<C_2$ such that, on the event $E$,
\begin{eqnarray*}
&&\Ecal(\ell\circ\hat{f})+
C_1 \Biggl(\sum_{j\notin J_f}\tau\bar{\epsilon}_j
\|\hat{f}_j\|_{L_2(\Pi)}
+\sum_{j=1}^N\tau^2\bar{\epsilon}_j^2\|\hat{f}_j\|_{\Hcal
_j} \Biggr)\\
&&\qquad\le\Ecal(\ell\circ f)
+C_2\sum_{j\in J_f} \bigl(\tau\bar{\epsilon}_j\|f_j-\hat{f}_j\|
_{L_2(\Pi)}
+\tau^2\bar{\epsilon}_j^2\|f_j\|_{\Hcal_j} \bigr)\\
&&\qquad\quad{} + |(P_n-P) (\ell\circ\hat{f}-\ell\circ f) |.
\end{eqnarray*}

We now bound the empirical process
$ |(P_n-P) (\ell\circ\hat{f}-\ell\circ f ) |$,
where we use the following result that will be proved in the next
section. Suppose that $f=\sum_{j=1}^N f_j$, $f_j\in\mathcal{H}_j$ and
$\|f\|_{L_{\infty}}\leq R$ (we will need it with $R=R_D^{\ast}$).
Denote
\begin{eqnarray*}
\Gcal(\Delta_-,\Delta_+, R)&=& \Biggl\{g\dvtx
\sum_{j=1}^N\bar{\epsilon}_j\|g_j-f_j\|_{L_2(\Pi)}\le\Delta_-, \\
&&\hspace*{6.9pt}\sum_{j=1}^N\bar{\epsilon}_j^2\|g_j-f_j\|_{\Hcal
_j}\le
\Delta_+, \Biggl\|\sum_{j=1}^Ng_j\Biggr\|_{L_{\infty}}\le R \Biggr\}.
\end{eqnarray*}

\begin{lemma}
\label{th:uniep}
There exists a numerical constant $C>0$
such that for an arbitrary $A\geq1$ involved
in the definition of $\bar\epsilon_j, j=1,\ldots, N$
with probability at least $1-2N^{-A/2}$,
for all
%
%
\begin{equation}
\Delta_-\leq e^N,\qquad \Delta_+\le e^N,
\end{equation}
the following bound holds:
%
%
\begin{equation}
{\sup_{g\in\Gcal(\Delta_-,\Delta_+, R_D^{\ast})}}
|(P_n-P) (\ell
\circ
g-\ell\circ f ) |\le CL_{\ast} (\Delta_-+\Delta_+
+ e^{-N} ).
\end{equation}
\end{lemma}

Assuming that
%
%
\begin{equation}
\label{case_1}
\sum_{j=1}^N\bar{\epsilon}_j\|\hat f_j-f_j\|_{L_2(\Pi)}\le e^N,\qquad
\sum_{j=1}^N\bar{\epsilon}_j^2\|\hat f_j-f_j\|_{\Hcal_j}\le e^N
\end{equation}
and using the lemma, we get
\begin{eqnarray*}
&&\Ecal(\ell\circ\hat{f})+
C_1 \Biggl(\sum_{j\notin J_f}\tau\bar{\epsilon}_j\|\hat{f}_j\|
_{L_2(\Pi)}
+\sum_{j=1}^N\tau^2\bar{\epsilon}_j^2\|\hat{f}_j\|_{\Hcal
_j} \Biggr)\\
&&\qquad\le\Ecal(\ell\circ f)+
C_2\sum_{j\in J_f} \bigl(\tau\bar{\epsilon}_j\|f_j-\hat{f}_j\|
_{L_2(\Pi)}
+\tau^2\bar{\epsilon}_j^2\|f_j\|_{\Hcal_j} \bigr)\\
&&\qquad\quad{}+C_3L_{\ast}\sum_{j=1}^N \bigl(\bar{\epsilon}_j\|\hat{f}_j-f_j\|
_{L_2(\Pi)}
+\bar{\epsilon}_j^2\|\hat{f}_j-f_j\|_{\Hcal_j} \bigr)+C_3L_{\ast
}e^{-N}\\
&&\qquad\le\Ecal(\ell\circ f)+
C_2\sum_{j\in J_f} \bigl(\tau\bar{\epsilon}_j\|f_j-\hat{f}_j\|
_{L_2(\Pi)}
+\tau^2\bar{\epsilon}_j^2\|f_j\|_{\Hcal_j} \bigr)\\
&&\qquad\quad{}+C_3L_{\ast}\sum_{j=1}^N \bigl(\bar{\epsilon}_j\|\hat{f}_j-f_j\|
_{L_2(\Pi)}
+\bar{\epsilon}_j^2\|\hat{f}_j\|_{\Hcal_j}
+\bar{\epsilon}_j^2\|f_j\|_{\Hcal_j} \bigr)\\
&&\qquad\quad{} + C_3L_{\ast}e^{-N}
\end{eqnarray*}
for some numerical constant $C_3>0$.
By choosing a numerical constant $B$ properly, $\tau$ can be made large
enough so
that $2C_3L_{\ast}\le\tau C_1\leq\tau C_2$. Then, we have
%
%
\begin{eqnarray}
\label{bound_X}
&&
\Ecal(\ell\circ\hat{f})+{1\over2}C_1 \Biggl(\sum_{j\notin J_f}
\tau\bar{\epsilon}_j\|\hat{f}_j\|_{L_2(\Pi)}
+\sum_{j=1}^N\tau^2\bar{\epsilon}_j^2\|\hat{f}_j\|_{\Hcal
_j} \Biggr)
\nonumber\\
&&\qquad\le\Ecal(\ell\circ f)+
2C_2\sum_{j\in J_f} \bigl(\tau\bar{\epsilon}_j\|f_j-\hat{f}_j\|
_{L_2(\Pi)}
+\tau^2\bar{\epsilon}_j^2\|f_j\|_{\Hcal_j} \bigr)\\
&&\qquad\quad{}+
(C_2/2) \tau e^{-N},\nonumber
\end{eqnarray}
which also implies
%
%
\begin{eqnarray}
\label{bound_Y}
&&\Ecal(\ell\circ\hat{f})+{1\over2}C_1 \Biggl(\sum_{j=1}^N
\tau\bar{\epsilon}_j\|\hat{f}_j-f_j\|_{L_2(\Pi)}
+\sum_{j=1}^N\tau^2\bar{\epsilon}_j^2\|\hat{f}_j\|_{\Hcal
_j} \Biggr)
\nonumber
\\
&&\qquad\le\Ecal(\ell\circ f)+
\biggl(2C_2+\frac{C_1}{2} \biggr)\sum_{j\in
J_f}\tau\bar{\epsilon}_j\|f_j-\hat{f}_j\|_{L_2(\Pi)}
\\
&&\qquad\quad{}+2C_2 \tau^2 \sum_{j\in J_f}\bar{\epsilon}_j^2\|f_j\|_{\Hcal_j}+
(C_2/2) \tau e^{-N}.\nonumber
\end{eqnarray}
We first consider the case when
%
%
\begin{eqnarray}
\label{bound_Z}
4 C_2\sum_{j\in J_f}\tau\bar{\epsilon}_j\|f_j-\hat{f}_j\|_{L_2(\Pi)}
&\ge&
\Ecal(\ell\circ f)+
2C_2\sum_{j\in
J_f}\tau^2\bar{\epsilon}_j^2\|f_j\|_{\Hcal_j}\nonumber\\[-8pt]\\[-8pt]
&&{}+(C_2/2)\tau e^{-N}.\nonumber
\end{eqnarray}
Then (\ref{bound_X}) implies that
%
%
\begin{eqnarray}
\label{eq:riskbddom}
&&\Ecal(\ell\circ\hat{f})+{1\over2}
C_1 \Biggl(\sum_{j\notin J_f}\tau\bar{\epsilon}_j\|\hat{f}_j\|
_{L_2(\Pi)}
+\sum_{j=1}^N\tau^2\bar{\epsilon}_j^2\|\hat{f}_j\|_{\Hcal
_j} \Biggr)\nonumber\\[-8pt]\\[-8pt]
&&\qquad\le6 C_2\sum_{j\in J_f}\tau\bar{\epsilon}_j\|f_j-\hat{f}_j\|
_{L_2(\Pi)},\nonumber
\end{eqnarray}
which yields
%
%
\begin{equation}
\sum_{j\notin J_f}\tau\bar{\epsilon}_j\|\hat{f}_j\|_{L_2(\Pi)}
\le{12 C_2\over C_1}\sum_{j\in J_f}\tau\bar{\epsilon}_j
\|f_j-\hat{f}_j\|_{L_2(\Pi)}.
\end{equation}
Therefore, $(\hat f_1-f_1,\ldots, \hat f_N-f_N)\in\mathcal{K}_{J_f}^{(b)}$
with $b:=12 C_2/C_1$.
Using the definition of $\beta_b(J_f)$, it follows
from (\ref{bound_Y}), (\ref{bound_Z}) and the assumption $C_1<1<C_2$ that
\begin{eqnarray*}
&&\Ecal(\ell\circ\hat{f})
+{1\over2}C_1 \Biggl(\sum_{j=1}^N\tau\bar{\epsilon}_j
\|\hat{f}_j-f_j\|_{L_2(\Pi)}
+\sum_{j=1}^N\tau^2\bar{\epsilon}_j^2\|\hat{f}_j\|_{\Hcal
_j} \Biggr)\\
&&\qquad\le \biggl(6C_2 +\frac{C_1}{2} \biggr)
\tau\beta_b(J_f)\|f-\hat{f}\|_{L_2(\Pi)}\\
&&\qquad\le 7 C_2\tau\beta_b(J_f) \bigl(\|f-f_\ast\|_{L_2(\Pi)}
+\|f_\ast-\hat{f}\|_{L_2(\Pi)} \bigr).
\end{eqnarray*}

Recall that for losses of quadratic type
%
%
\begin{equation}
\label{quadr}\quad
\Ecal(\ell\circ f)\ge m_{\ast} \|f-f_\ast\|_{L_2(\Pi)}^2
\quad\mbox{and}\quad
\Ecal(\ell\circ\hat{f})\ge m_{\ast} \|\hat{f}-f_\ast\|_{L_2(\Pi)}^2.
\end{equation}
%
Then
\begin{eqnarray*}
&&\Ecal(\ell\circ\hat{f})
+{1\over2}C_1 \Biggl(\sum_{j=1}^N \tau\bar{\epsilon}_j
\|\hat{f}_j-f_j\|_{L_2(\Pi)}
+\sum_{j=1}^N\tau^2\bar{\epsilon}_j^2\|\hat{f}_j\|_{\Hcal
_j} \Biggr)\\
&&\qquad\le 7 \tau C_2 m_{\ast}^{-1/2}\beta_b(J_f) \bigl(\Ecal^{1/2}(\ell
\circ
f)+\Ecal^{1/2}(\ell\circ\hat{f}) \bigr).
\end{eqnarray*}
Using the fact that $ab\le(a^2+b^2)/2$, we get
%
%
\begin{equation}\qquad
7\tau C_2m_{\ast}^{-1/2}\beta_b(J_f)\Ecal^{1/2}(\ell\circ f)\le
(49/2)\tau^2C_2^2m_{\ast}^{-1}\beta_b^2(J_f)+\tfrac{1}{2}\Ecal(\ell
\circ f)
\end{equation}
and
%
%
\begin{equation}\qquad
7\tau C_2 m_{\ast}^{-1/2}\beta_b(J_f)\Ecal^{1/2}(\ell\circ\hat
{f})\le
(49/2)\tau^2C_2^2m_{\ast}^{-1}\beta_b^2(J_f)+\tfrac{1}{2}\Ecal(\ell
\circ
\hat{f}).
\end{equation}
Therefore,
%
%
\begin{eqnarray}
\label{eq:riskbddom1}
&&\Ecal(\ell\circ\hat{f})+C_1\sum_{j=1}^N\tau
\bar{\epsilon}_j\|\hat{f}_j\|_{L_2(\Pi)}+C_1\sum_{j=1}^N\tau
^2\bar
{\epsilon}_j^2\|\hat{f}_j\|_{\Hcal_j}\nonumber\\[-8pt]\\[-8pt]
&&\qquad
\le\Ecal(\ell\circ f)+100
\tau^2C_2^2m_{\ast}^{-1}\beta_b^2(J_f).\nonumber
\end{eqnarray}

We now consider the case when
%
%
\begin{eqnarray}
&&
4 C_2\sum_{j\in J_f}\tau\bar{\epsilon}_j\|f_j-\hat{f}_j\|_{L_2(\Pi)}
\nonumber\\[-8pt]\\[-8pt]
&&\qquad< \Ecal(\ell\circ f)+
2C_2\sum_{j\in J_f}\tau^2\bar{\epsilon}_j^2\|f_j\|_{\Hcal
_j}+(C_2/2)\tau e^{-N}.\nonumber
\end{eqnarray}
It is easy to derive from (\ref{bound_Y}) that in this case
%
%
\begin{eqnarray}\qquad
\label{eq:riskbddom2}
&&
\Ecal(\ell\circ\hat{f})+
{1\over2}C_1 \Biggl(\sum_{j=1}^N
\tau\bar{\epsilon}_j\|\hat{f}_j-f_j\|_{L_2(\Pi)}
+\sum_{j=1}^N\tau^2\bar{\epsilon}_j^2\|\hat{f}_j\|_{\Hcal
_j} \Biggr)
\nonumber\\[-8pt]\\[-8pt]
&&\qquad
\le
\biggl(\frac{3}{2}+\frac{C_1}{8C_2} \biggr) \biggl(\Ecal(\ell\circ f)+
2C_2\sum_{j\in J_f}\tau^2\bar{\epsilon}_j^2\|f_j\|_{\Hcal
_j}+(C_2/2)\tau
e^{-N} \biggr).\nonumber
\end{eqnarray}
Since $\beta_b(J_f)\geq\sqrt{\frac{A\log N}{n}}$ [see the comment after
the definition
of $\beta_b(J_f)$],
we have
\[
\tau e^{-N}
\leq\tau^2 \sqrt{\frac{A\log N}{n}}\leq\tau^2 \beta_b^2(J_f),
\]
where we also used the assumptions that $\log N \geq2\log\log n$ and
$A\geq4$.
Substituting this in (\ref{eq:riskbddom2}) and then combining the
resulting bound with (\ref{eq:riskbddom1}) concludes
the proof of (\ref{main_oracle}) in the case when conditions (\ref{case_1})
hold.

It remains to consider the case when (\ref{case_1}) does not hold.
The main idea is to show that in this case the right-hand side
of the oracle inequality is rather large while we still can control the
left-hand
side, so, the inequality becomes
trivial.
To this end, note that, by the definition of $\hat f$,
for some numerical constant $c_1$,
\[
P_n(\ell\circ\hat f)
+\sum_{j=1}^N \bigl(\tau\hat{\epsilon}_j\|\hat{f}_j\|_{L_2(\Pi_n)}
+\tau^2\hat{\epsilon}_j^2\|\hat{f}_j\|_{\Hcal_j} \bigr)
\leq n^{-1} \sum_{j=1}^n\ell(Y_j;0)\leq c_1
\]
[since the value of the penalized empirical risk at $\hat f$ is
not larger than its value at $f=0$ and,
by the assumptions on the loss, $\ell(y,0)$ is
uniformly bounded by a numerical constant].
The last equation implies
that, on the event $E$ defined
earlier in the proof
[see (\ref{condition_01}), (\ref{condition_02})], the following bound holds:
\[
\sum_{j=1}^N \frac{\tau}{C}\bar{\epsilon}_j \biggl(\frac{1}{C}\|
\hat
{f}_j\|_{L_2(\Pi)}-
\bar\epsilon_j \|\hat f_j\|_{\mathcal{H}_j} \biggr)
+\sum_{j=1}^N \frac{\tau^2}{C^2}\bar{\epsilon}_j^2\|\hat{f}_j\|
_{\Hcal_j}
\leq c_1.
\]
Equivalently,
\[
\frac{\tau}{C^2}\sum_{j=1}^N \bar{\epsilon}_j\|\hat{f}_j\|
_{L_2(\Pi)}
+
\biggl(\frac{\tau^2}{C^2}-\frac{\tau}{C} \biggr)
\sum_{j=1}^N \bar{\epsilon}_j^2\|\hat{f}_j\|_{\Hcal_j}
\leq c_1.
\]
As soon as $\tau\geq2C$, so that $\tau^2/C^2-\tau/C\geq\tau^2/(2C^2)$,
we have
%
%
\begin{equation}
\label{hatf}
\tau\sum_{j=1}^N \bar{\epsilon}_j\|\hat{f}_j\|_{L_2(\Pi)}
+
\tau^2\sum_{j=1}^N \bar{\epsilon}_j^2\|\hat{f}_j\|_{\Hcal_j}
\leq2c_1 C^2.
\end{equation}
Note also that, by the assumptions on the loss function,
%
%
\begin{eqnarray}
\label{excess_A}
\mathcal{E}(\ell\circ\hat f)&\leq&
P(\ell\circ\hat f)\nonumber\\
&\leq&
{\mathbb E}\ell(Y;0)+
|P(\ell\circ\hat f)-P(\ell\circ0)|\nonumber\\[-8pt]\\[-8pt]
&\leq& c_1 + L_{\ast}\|\hat f\|_{L_{2}(\Pi)} \leq
c_1 +L_{\ast} \sum_{j=1}^N \|\hat f\|_{L_2(\Pi)}
\nonumber\\
&\leq&
c_1 + 2c_1 C^2 L_{\ast} \frac{1}{\tau}\sqrt{\frac{n}{A\log
N}},\nonumber
\end{eqnarray}
where we used the Lipschitz condition on $\ell$, and
also bound (\ref{hatf}) and the fact that
$\bar\epsilon_j\geq\sqrt{A\log N/n}$ (by its definition).

Recall that we are considering the case when (\ref{case_1}) does not
hold. We will consider two cases: (a) when $e^N\leq c_3$, where
$c_3\geq c_1$
is a numerical constant, and (b)
when $e^N>c_3$. The first case is very simple since $N$ and $n$ are
both upper bounded by a numerical constant
(recall the assumption $\log N\geq2\log\log n$). In this case,
$\beta_b(J_f)\geq\sqrt{\frac{A\log N}{n}}$
is bounded from below\vspace*{2pt} by a numerical constant.
As a consequence of these observations, bounds~(\ref{hatf}) and (\ref
{excess_A})
imply that\vspace*{-2pt}
\begin{eqnarray*}
&&
\Ecal(\ell\circ\hat{f})+C_1 \Biggl(\sum_{j=1}^N
\tau\bar{\epsilon}_j\|\hat{f}_j\|_{L_2(\Pi)}
+\sum_{j=1}^N\tau^2\bar{\epsilon}_j^2\|\hat{f}_j\|_{\Hcal
_j} \Biggr)\\[-3pt]
&&\qquad\leq C_2 \tau^2 \beta_b^2(J_f)\vspace*{-2pt}
\end{eqnarray*}
for some numerical constant $C_2>0$. In the case (b), we have\vspace*{-2pt}
\[
\sum_{j=1}^N \bar\epsilon_j \|\hat f_j-f_j\|_{L_2(\Pi)}+
\sum_{j=1}^N \bar\epsilon_j^2 \|\hat f_j-f_j\|_{\mathcal{H}_j}\geq e^N\vspace*{-2pt}
\]
and, in view of (\ref{hatf}), this implies\vspace*{-2pt}
\[
\sum_{j=1}^N \bar\epsilon_j \|f_j\|_{L_2(\Pi)}+
\sum_{j=1}^N \bar\epsilon_j^2 \|f_j\|_{\mathcal{H}_j}
\geq e^N -c_1/2\geq e^N/2.\vspace*{-2pt}
\]
So, either we have\vspace*{-2pt}
\[
\sum_{j=1}^N \bar\epsilon_j^2 \|f_j\|_{\mathcal{H}_j}\geq e^N/4\vspace*{-2pt}
\]
or\vspace*{-2pt}
\[
\sum_{j=1}^N \bar\epsilon_j \|f_j\|_{L_2(\Pi)}\geq e^N/4.\vspace*{-2pt}
\]
Moreover, in the second case, we also have\vspace*{-2pt}
\begin{eqnarray*}
\sum_{j=1}^N \bar\epsilon_j^2 \|f_j\|_{\mathcal{H}_j}&\geq&
\sqrt{\frac{A\log N}{n}}
\sum_{j=1}^N \bar\epsilon_j \|f_j\|_{L_2(\Pi)}\\[-3pt]
&\geq&(e^N/4)\sqrt
{\frac
{A\log N}{n}}.\vspace*{-2pt}
\end{eqnarray*}
In both cases we can conclude that, under the assumption that $\log
N\geq2\log\log n$ and $e^N>c_3$
for a sufficiently large numerical constant $c_3$,
\begin{eqnarray*}
&&\mathcal{E}(\ell\circ\hat f)
+\sum_{j=1}^N \bigl(\tau\bar{\epsilon}_j\|\hat{f}_j\|_{L_2(\Pi)}
+\tau^2\bar{\epsilon}_j^2\|\hat{f}_j\|_{\Hcal_j} \bigr)
\\[-3pt]
&&\qquad\leq
c_1 + 2c_1 C^2 L_{\ast} \frac{1}{\tau}\sqrt{\frac{n}{A\log N}}+
2c_1 C^2\\
&&\qquad\leq\frac{\tau^2 e^N}{4}\sqrt{\frac{A\log N}{n}}\leq\tau^2
\sum_{j\in J_f} \bar\epsilon_j^2 \|f_j\|_{\mathcal{H}_j}.
\end{eqnarray*}
Thus, in both cases (a) and (b), the following bound holds:
%
%
\begin{eqnarray}
\label{konec}
&&\Ecal(\ell\circ\hat{f})+C_1 \Biggl(\sum_{j=1}^N
\tau\bar{\epsilon}_j\|\hat{f}_j\|_{L_2(\Pi)}
+\sum_{j=1}^N\tau^2\bar{\epsilon}_j^2\|\hat{f}_j\|_{\Hcal
_j} \Biggr)\nonumber\\[-8pt]\\[-8pt]
&&\qquad
\leq
C_2 \tau^2 \biggl(\sum_{j\in J_f} \bar\epsilon_j^2 \|f_j\|
_{\mathcal
{H}_j}+\beta_b^2(J_f)
\biggr).\nonumber
\end{eqnarray}

To complete the proof, observe that
%
%
\begin{eqnarray}
\label{konec1}
&&
\Ecal(\ell\circ\hat{f})+C_1 \Biggl(\sum_{j=1}^N
\tau\bar{\epsilon}_j\|\hat{f}_j-f_j\|_{L_2(\Pi)}
+
\sum_{j=1}^N\tau^2\bar{\epsilon}_j^2\|\hat{f}_j\|_{\Hcal_j} \Biggr)
\nonumber
\\
&&\qquad
\leq
\Ecal(\ell\circ\hat{f})+
C_1 \Biggl(\sum_{j=1}^N
\tau\bar{\epsilon}_j\|\hat{f}_j\|_{L_2(\Pi)}
+\sum_{j=1}^N\tau^2\bar{\epsilon}_j^2\|\hat{f}_j\|_{\Hcal
_j} \Biggr)\nonumber\\
&&\qquad\quad{} +C_1 \sum_{j\in J_f}\tau\bar\varepsilon_j \|\hat f_j-f_j\|_{L_2(\Pi)}
\\
&&\qquad
\leq C_2 \tau^2 \biggl(\sum_{j\in J_f} \bar\epsilon_j^2 \|f_j\|
_{\mathcal{H}_j}+\beta_b^2(J_f)
\biggr)\nonumber\\
&&\qquad\quad{}
+C_2 \sum_{j\in J_f}\tau\bar\varepsilon_j \|\hat
f_j-f_j\|_{L_2(\Pi)}.\nonumber
\end{eqnarray}
Note also that, by the definition of $\beta_b(J_f)$, for all
$b>0$,
%
%
\begin{eqnarray}\quad
&&
\sum_{j\in J_f}\tau\bar\varepsilon_j \|\hat f_j-f_j\|_{L_2(\Pi)}
\nonumber\\
&&\qquad\leq
\tau\beta_b(J_f) \biggl\|\sum_{j\in J_f}(\hat f_j-f_j) \biggr\|
_{L_2(\Pi
)}\nonumber\\[-8pt]\\[-8pt]
&&\qquad\leq
\tau\beta_b(J_f)\|\hat f-f\|_{L_2(\Pi)}+
\tau\beta_b(J_f)\sqrt{\frac{n}{A\log N}}\sum_{j\notin J_f}\bar
\varepsilon_j \|\hat f_j\|_{L_2(\Pi)}
\nonumber\\
&&\qquad\leq
\tau\beta_b(J_f)\|\hat f-f\|_{L_2(\Pi)}+
\tau\beta_b(J_f)
\frac{2c_1C^2}{\tau}
\sqrt{\frac{n}{A\log N}},\nonumber
\end{eqnarray}
where we used the fact that, for all $j$, $\bar\varepsilon_j \geq
\sqrt
{\frac{A\log N}{n}}$
and also bound (\ref{hatf}).
By an argument similar
to (\ref{quadr})--(\ref{eq:riskbddom1}), it is easy to deduce from the
last bound that
%
%
\begin{eqnarray}\qquad
C_2\sum_{j\in J_f}\tau\bar\varepsilon_j \|\hat f_j-f_j\|_{L_2(\Pi)}
&\leq&\frac{3}{2}\frac{C_2^2 \tau^2}{m_{\ast}} \beta_b^2(J_f)
+\frac{1}{2}\mathcal{E}(\ell\circ\hat f)
+\frac{1}{2}\mathcal{E}(\ell\circ f)\nonumber\\[-9pt]\\[-9pt]
&&{}+\frac{2c_1^2C^4}{\tau^2}\frac{n}{A\log N}.\nonumber
\end{eqnarray}
Substituting this in bound (\ref{konec1}), we get
%
%
\begin{eqnarray}
&&
\frac{1}{2}\Ecal(\ell\circ\hat{f})+C_1 \Biggl(\sum_{j=1}^N
\tau\bar{\epsilon}_j\|\hat{f}_j-f_j\|_{L_2(\Pi)}
+
\sum_{j=1}^N\tau^2\bar{\epsilon}_j^2\|\hat{f}_j\|_{\Hcal_j} \Biggr)
\nonumber\\[-2pt]
&&\qquad
\leq
C_2 \tau^2 \biggl(\sum_{j\in J_f} \bar\epsilon_j^2 \|f_j\|
_{\mathcal{H}_j}+
\beta_b^2(J_f)
\biggr)
\nonumber\\[-2pt]
&&\qquad\quad{}
+\frac{3}{2}\frac{C_2^2 \tau^2}{m_{\ast}} \beta_b^2(J_f)
+\frac{1}{2}\mathcal{E}(\ell\circ f)+
\frac{2c_1^2C^4}{\tau^2}\frac{n}{A\log N}
\\[-2pt]
&&\qquad
\leq
\frac{1}{2}\mathcal{E}(\ell\circ f)+
C_2^{\prime} \tau^2 \biggl(\sum_{j\in J_f} \bar\epsilon_j^2 \|
f_j\|
_{\mathcal{H}_j}+
\frac{\beta_b^2(J_f)}{m_{\ast}}
\biggr)\nonumber\\[-2pt]
&&\qquad\quad{}+\frac{2c_1^2C^2}{\tau^2}\frac{n}{A\log N},\nonumber
\end{eqnarray}
with some numerical constant $C_2^{\prime}$.
It is enough now to observe [considering again the cases (a) and (b),
as it was done before], that either the last term is upper bounded by
$\sum_{j\in J_f}\bar\varepsilon_j\|f_j\|_{\mathcal{H}_j}$, or it is upper
bounded by $\beta_b^2(J_f)$, to complete the proof.
\end{pf}

Now, to derive Theorem~\ref{co:homobd}, it is enough to check
that, for a numerical constant $c>0$,
\begin{eqnarray*}
\beta_b(J_f)&\leq&
\biggl(\sum_{j\in J_{f}}\bar\epsilon_j^2 \biggr)^{1/2}
\beta_{2,\infty}(J_f)
\\[-2pt]
&\leq&
c \biggl(\sum_{j\in J_{f}}\breve\epsilon_j^2 \biggr)^{1/2}
\beta_{2,\infty}(J_f),
\end{eqnarray*}
which easily follows from the definitions of $\beta_b$ and $\beta
_{2,\infty}$.
Similarly, the proof of
Theorem~\ref{co:homobd-1} follows from the fact that, under the assumption
that
$
\Lambda^{-1}\leq\frac{\breve\epsilon_j}{\breve\epsilon}\leq
\Lambda,
$
we have
$
\mathcal{K}_J^{(b)}\subset K_J^{(b^{\prime})},
$
where $b^{\prime}=c \Lambda^2 b$, $c$ being a numerical constant.
This easily implies the bound
$
\beta_b (J_f)\leq c_1\Lambda\beta_{2,b^{\prime}}(J_f)
\sqrt{d(f)} \breve\epsilon,
$
where $c_1$ is a numerical constant.\vadjust{\goodbreak}

\section{Bounding the empirical process}
\label{sec:empirical_process}

We now proceed to prove
Lemma~\ref{th:uniep} that was used
to bound
$ |(P_n-P) (\ell\circ\hat{f}-\ell\circ f ) |$.
To this end, we begin with a fixed pair
$(\Delta_-, \Delta_+)$. Throughout the proof, we write
$R:=R_D^{\ast}$.
By Talagrand's concentration inequality, with probability at least $1-e^{-t}$
\begin{eqnarray*}
&&{\sup_{g\in\Gcal(\Delta_-, \Delta_+, R)} }|(P_n-P)
(\ell\circ g-\ell\circ f ) |
\\[-3pt]
&&\qquad\le
2 \Biggl(\bE\Bigl[{\sup_{g\in\Gcal(\Delta_-, \Delta_+,R)}}
|(P_n-P) (\ell\circ g-\ell\circ f ) | \Bigr]\\[-3pt]
&&\qquad\quad\hspace*{9.3pt}{}+\|\ell\circ g-\ell\circ f\|_{L_2(P)}
\sqrt{t\over n}+\|\ell\circ g-\ell\circ f\|_{L_{\infty}}{t\over
n} \Biggr).
\end{eqnarray*}
Now note that
\begin{eqnarray*}
\|\ell\circ g-\ell\circ f\|_{L_2(P)}
&\le&
L_{\ast}\|g-f\|_{L_2(\Pi)}
\\[-3pt]
&\le& L_{\ast}\sum_{j=1}^N\|g_j-f_j\|_{L_2(\Pi)}\\[-3pt]
&\le&
L_{\ast} \Bigl(\min_j\bar{\epsilon}_j \Bigr)^{-1}
\sum_{j=1}^N\bar{\epsilon}_j\|g_j-f_j\|_{L_2(\Pi)},
\end{eqnarray*}
where we used the fact that
the Lipschitz
constant of the loss $\ell$ on the range of functions from
$\mathcal{G}(\Delta_-, \Delta_+, R)$ is bounded by $L_{\ast}$.
Together with the fact that
$\bar{\epsilon}_j\ge(A\log N/n)^{1/2}$ for all $j$, this yields
%
%
\begin{equation}
\|\ell\circ g-\ell\circ f\|_{L_2(P)}\le L_{\ast}\sqrt{\frac
{n}{A\log
N}}\Delta_-.
\end{equation}
Furthermore,
\begin{eqnarray*}
\|\ell\circ g-\ell\circ f\|_{L_{\infty}}
&\le& L_{\ast} \|g-f\|_{L_{\infty}}\\[-3pt]
&\le&
L_{\ast}\sum_{j=1}^N \|g_j-f_j\|_{\mathcal{H}_j}\\[-3pt]
&\le& L_{\ast} \frac{n}{A\log N}\Delta_{+}.
\end{eqnarray*}
%
In summary, we have
\begin{eqnarray*}
&&{\sup_{g\in\Gcal(\Delta_-, \Delta_+, R)}} |(P_n-P) (\ell
\circ
g-\ell\circ f ) |\\[-3pt]
&&\qquad\le
2 \Biggl(\bE\Bigl[{\sup_{g\in\Gcal(\Delta_-, \Delta_+,R)}}
|(P_n-P) (\ell\circ g-\ell\circ f ) |
\Bigr]\\[-3pt]
&&\qquad\quad\hspace*{27pt}{}+ L_{\ast}\Delta_-\sqrt{t\over A\log N}+L_{\ast}\Delta
_+{t\over
n}{n\over A\log N} \Biggr).
\end{eqnarray*}
Now, by symmetrization inequality,
%
%
\begin{eqnarray}
&&
\bE\Bigl[{\sup_{g\in\Gcal(\Delta_-, \Delta_+,R)}}
|(P_n-P) (\ell
\circ g-\ell\circ f ) | \Bigr]\nonumber\\[-8pt]\\[-8pt]
&&\qquad\le{2\bE\sup_{g\in\Gcal
(\Delta
_-, \Delta_+,R)}} |R_n(\ell\circ g-\ell\circ f) |.\nonumber
\end{eqnarray}
An application of Rademacher contraction inequality further yields
%
%
\begin{eqnarray}
&&\bE\Bigl[{\sup_{g\in\Gcal(\Delta_-, \Delta_+,R)}}
|(P_n-P) (\ell\circ g-\ell\circ f ) |
\Bigr]\nonumber\\[-8pt]\\[-8pt]
&&\qquad\le
CL_{\ast}\bE\sup_{g\in\Gcal(\Delta_-, \Delta_+,R)}
|R_n(g-f) |,\nonumber
\end{eqnarray}
where $C>0$ is a numerical constant [again, it was used here that
the Lipschitz
constant of the loss $\ell$ on the range of functions from
$\mathcal{G}(\Delta_-, \Delta_+, R)$ is bounded by $L_{\ast}$].
Applying Talagrand's concentration inequality another time, we get that
with probability at least $1-e^{-t}$
%
\begin{eqnarray*}
{\bE\sup_{g\in\Gcal(\Delta_-, \Delta_+,R)}} |R_n(g-f)|
&\le& C
\Biggl({\sup_{g\in\Gcal(\Delta_-, \Delta_+,R)}}
|R_n(g-f) |\\
&&\hspace*{12.3pt}{} +\Delta_-\sqrt{t\over A\log N}+
\Delta_+{t\over n}{n\over A\log N} \Biggr)
\end{eqnarray*}
for some numerical constant $C>0$.

Recalling the definition of $\check\epsilon_j:=\check\epsilon(K_j)$,
we get
%
%
\begin{equation}
|R_n(h_j)|\le\check{\epsilon}_j\|h_j\|_{L_2(\Pi)} +\check{\epsilon
}_j^2\|h_j\|_{\Hcal_j},\qquad
h_j\in\mathcal{H}_j.
\end{equation}
Hence, with probability at least $1-2e^{-t}$ and with some numerical
constant $C>0$
\begin{eqnarray*}
&&{\sup_{g\in\Gcal(\Delta_-, \Delta_+,R)}}
|(P_n-P) (\ell\circ g-\ell\circ f ) |\\
&&\qquad\le CL_{\ast} \Biggl({\sup_{g\in\Gcal(\Delta_-, \Delta_+,R)}}
|R_n(g-f) |+
\Delta_-\sqrt{t\over A\log N}+
\Delta_+{t\over n}{n\over A\log N} \Biggr)\\
&&\qquad\le CL_{\ast} \Biggl(\sup_{g\in\Gcal(\Delta_-, \Delta_+,R)}
\sum_{j=1}^N |R_n(g_j-f_j) |
+\Delta_-\sqrt{t\over A\log N}+
\Delta_+{t\over n}{n\over A\log N} \Biggr)\\
&&\qquad\le
CL_{\ast} \Biggl(\sup_{g\in\Gcal(\Delta_-, \Delta_+,R)}
\sum_{j=1}^N \bigl(\check{\epsilon}_j\|g_j-f_j\|_{L_2(\Pi)}
+\check{\epsilon}_j^2\|g_j-f_j\|_{\Hcal_j} \bigr)\\
&&\qquad\quad\hspace*{121.6pt}{} +
\Delta_-\sqrt{t\over A\log N}
+\Delta_+{t\over n}{n\over A\log N} \Biggr).
\end{eqnarray*}

Using (\ref{conclusion_B}), $\check\epsilon_j$ can be upper bounded
by $c\bar\epsilon_j$ with some numerical constant $c>0$ on an event $E$
of probability at least $1-N^{-A/2}$.
Therefore, the following bound is obtained:
\begin{eqnarray*}
&&{\sup_{g\in\Gcal(\Delta_-, \Delta_+,R)} }|(P_n-P)
(\ell\circ g-\ell\circ f ) |\\
&&\qquad\le CL_{\ast} \Biggl(\Delta_- +\Delta_++\Delta_-\sqrt{t\over
A\log N}+
\Delta_+{t\over n}{n\over A\log N} \Biggr).
\end{eqnarray*}
It holds on the event
$E\cap F(\Delta_-, \Delta_+, t)$,
where
${\mathbb P}(F(\Delta_-, \Delta_+, t))\geq1-2e^{-t}$.

We will now choose $t=A\log N+4\log N+4\log(2/\log2)$ and obtain a
bound that
holds uniformly over
%
%
\begin{equation}
\label{eq:bdDel}
e^{-N}\le\Delta_-\le e^N \quad\mbox{and}\quad e^{-N}\le\Delta_+\le e^N.
\end{equation}
To this end, consider
%
%
\begin{equation}
\Delta_j^{-}=\Delta_j^+:=2^{-j}.
\end{equation}
For any $\Delta_j^-$ and $\Delta_k^+$ satisfying (\ref{eq:bdDel}),
we have
\begin{eqnarray*}
&&{\sup_{g\in\Gcal(\Delta_j^-, \Delta_k^+,R)}}
|(P_n-P) (\ell\circ g-\ell\circ f ) |\\
&&\qquad\le CL_{\ast} \Biggl(\Delta_j^- +\Delta_k^+ +\Delta_j^-\sqrt
{t\over
A\log N}+\Delta_k^+\frac{t}{n}
{n\over A\log N} \Biggr)
\end{eqnarray*}
on the event
$E\cap F(\Delta_j^-, \Delta_k^+, t)$.
Therefore, simultaneously for all $\Delta_j^-$ and $\Delta_k^+$ satisfying
(\ref{eq:bdDel}),
we have
\begin{eqnarray*}
&&{\sup_{g\in\Gcal(\Delta_j^-, \Delta_k^+,R)}}
|(P_n-P) (\ell\circ g-\ell\circ f ) |\\
&&\qquad\le CL_{\ast} \Biggl(\Delta_j^- +\Delta_k^+
+\Delta_j^-\sqrt{A\log N+4\log N+4\log(2/\log2)\over A\log N}\\
&&\qquad\quad\hspace*{45.2pt}{}+\Delta_+{A\log N+4\log N+4\log(2/\log2)\over
n}{n\over
A\log N} \Biggr)
\end{eqnarray*}
on the event
$E^{\prime}:=E\cap(\bigcap_{j,k} F(\Delta_j^-, \Delta_k^+,
t) )$.
The last intersection is over all $j,k$ such that
conditions (\ref{eq:bdDel}) hold for $\Delta_j^-, \Delta_k^+$. The
number of the events
in this intersection is bounded by $(2/\log2)^2 N^2$. Therefore,
%
%
\begin{eqnarray}
{\mathbb P}(E^{\prime})&\geq&1-(2/\log2)^2 N^2\exp\bigl(-A\log
N-4\log
N-4\log(2/\log2) \bigr)\nonumber\\
&&{} - {\mathbb P}(E)\\
&\geq&1-2 N^{-A/2}.\nonumber
\end{eqnarray}
Using monotonicity of the functions of $\Delta_-, \Delta_+$
involved in the inequalities, the bounds can be extended
to the whole range of values of $\Delta_-, \Delta_+$ satisfying
(\ref{eq:bdDel}), so, with probability at least $1-2N^{-A/2}$ we
have for all such $\Delta_-, \Delta_+$
%
%
\begin{equation}
{\sup_{g\in\Gcal(\Delta_-,\Delta_+, R)} }|(P_n-P) (\ell
\circ
g-\ell\circ f ) |\le CL_{\ast} (\Delta_-+\Delta
_+ ).
\end{equation}
If $\Delta_- \leq e^{-N}$, or $\Delta_+\leq e^{-N}$,
it follows by monotonicity of the left-hand side that with the same
probability
%
%
\begin{equation}\quad
{\sup_{g\in\Gcal(\Delta_-,\Delta_+, R)}} |(P_n-P) (\ell
\circ
g-\ell\circ f ) |\le CL_{\ast} (\Delta_-+\Delta_+
+e^{-N} ),
\end{equation}
which completes the proof.

\section*{Acknowledgments}
The authors are thankful to the referees for a number of helpful
suggestions. The first author is thankful to Evarist Gin\'e for useful
conversations about the paper.


%
\printaddresses

\end{document}